\documentclass{jT}

\usepackage{amssymb,amsmath,latexsym,amsthm}
\usepackage[english]{babel}

\newtheorem{theorem}{Theorem}[section]
\newtheorem{lemma}{Lemma}[section]
\newtheorem{proposition}{Proposition}[section]

\theoremstyle{definition}
\newtheorem*{definition}{Definition}
\newtheorem*{remark}{Remark}

\numberwithin{equation}{section}

\renewcommand{\leq}{\leqslant}
\renewcommand{\geq}{\geqslant}

\newsavebox{\proofbox}
\savebox{\proofbox}{\begin{picture}(7,7)%
  \put(0,0){\framebox(7,7){}}\end{picture}}

\newcommand{\md}[1]{\ensuremath{(\mbox{mod}\, #1)}}
\newcommand{\mdsub}[1]{\ensuremath{(\mbox{\scriptsize mod}\, #1)}}

\newcommand{\mdsublem}[1]{\ensuremath{(\mbox{\scriptsize \textup{mod}}\, #1)}}

\newcommand{\lff}{{\fontencoding{T1}\fontfamily{ptm}\selectfont<<}}
\newcommand{\rff}{{\fontencoding{T1}\fontfamily{ptm}\selectfont>>}}

\def\endproof{\hfill{\usebox{\proofbox}}}

\def\E{\mathbb{E}}
\def\Z{\mathbb{Z}}
\def\N{\mathbb{N}}
\def\R{\mathbb{R}}
\def\C{\mathbb{C}}
\def\lamsel{\lambda^{\mbox{\scriptsize SEL}}}
\def\lamgy{\lambda^{\mbox{\scriptsize GY}}}
\def\ni{\noindent}
\def\vs{\vspace{11pt}}

\begin{document}

\title[Restriction theory of the Selberg sieve]{Restriction theory of the Selberg sieve, with applications}


\author[Ben {\sc Green}]{{\sc Ben} GREEN}
\address{Ben {\sc Green}\\
School of Mathematics\\
University of Bristol\\
Bristol BS8 1TW\\
England}
\email{b.j.green@bristol.ac.uk}
\urladdr{http://www.maths.bris.ac.uk/$\widetilde{\;}$mabjg}

\author[Terence {\sc Tao}]{{\sc Terence} TAO}
\address{Terence {\sc Tao}\\
Department of Mathematics\\
University of California at Los Angeles\\ 
Los Angeles CA 90095}
\email{tao@math.ucla.edu}
\urladdr{http://www.math.ucla.edu/$\widetilde{\;}$tao}

\maketitle

\begin{resume}
\ni Le crible de Selberg fournit des majorants pour certaines suites arithm\'etiques, comme les nombres premiers et les nombres premiers jumeaux. Nous demontrons un th\'eor\`eme de restriction $L^2$-$L^p$ pour les majorants de ce type. Comme application imm\'ediate, nous consid\'erons l'estimation des sommes exponentielles sur les $k$-uplets premiers. Soient les entiers positifs $a_1,\dots,a_k$ et $b_1,\dots,b_k$. Ecrit $h(\theta) := \sum_{n \in X} e(n\theta)$, ou $X$ est l'ensemble de tous $n \leq N$ tel que tous les nombres $a_1n + b_1,\dots, a_kn + b_k$ sont premiers. Nous obtenons les bornes sup\'erieures pour $\Vert h \Vert_{L^p(\mathbb{T})}$, $p > 2$, qui sont (en supposant la v\'erit\'e de la conjecture de Hardy et Littlewood sur les $k$-uplets premiers) d'ordre de magnitude correct. Une autre application est la suivante. En utilisant les th\'eor\`emes de Chen et de Roth et un \lff principe de transf\'erence \rff, nous demontrons qu'il existe une infinit\'e de suites arithm\'etiques $p_1 < p_2 < p_3$ de nombres premiers, telles que chacun $p_i + 2$ est premier ou un produit de deux nombres premier.
\end{resume}

\begin{abstr}
The Selberg sieve provides majorants for certain arithmetic sequences, such as the primes and the twin primes. 
We prove an $L^2$--$L^p$ restriction theorem for majorants of this type. An immediate application is to the estimation of exponential sums over prime $k$-tuples. Let $a_1,\dots,a_k$ and $b_1,\dots,b_k$ be positive integers. Write $h(\theta) := \sum_{n \in X} e(n\theta)$, where $X$ is the set of all $n \leq N$ such that the numbers $a_1n + b_1,\dots, a_kn + b_k$ are all prime. We obtain upper bounds for $\Vert h \Vert_{L^p(\mathbb{T})}$, $p > 2$, which are (conditionally on the Hardy-Littlewood prime tuple conjecture) of the correct order of magnitude. As a second application we deduce from Chen's theorem, Roth's theorem, and a transference principle that there are infinitely many arithmetic progressions $p_1 < p_2 < p_3$ of primes, such that $p_i + 2$ is either a prime or a product of two primes for each $i=1,2,3$.
\end{abstr}

\bigskip

\section{Introduction and statement of results} 
\ni The exponential sum over the primes\footnote{The letter $p$ will be used both for primes and for exponents in H\"older's inequality; there should be no danger of confusion.  All summations over $p$ are over the primes, and summations over other variables are over the positive integers unless otherwise indicated. Most of the notation in the introduction is very standard -- in the event of further clarification being required, consultation of \S \ref{notation-sec} is recommended.}, 
\[ h_N(\theta)  :=  \sum_{p \leq N} e(p\theta)\] 
has been very extensively studied.  Here of course $e(\theta) := e^{2\pi i \theta}$.  Since the time of Hardy and Littlewood such sums have been treated by dividing the circle $\R/\Z$ into the \textit{major arcs}, which is to say points near a rational with small denominator, and the complement of this set, the \textit{minor arcs}.\vs

\ni The sum $h_N(\theta)$ (or at least the very closely-related exponential sum over the von Mangoldt function) is discussed in \cite[Ch. 3]{Vaughan}, for example. There it is shown how to use information concerning the distribution of primes in arithmetic progressions in order to get an asymptotic for $h_N(\theta)$ on the major arcs. On the minor arcs, $h_N(\theta)$ is shown to be rather small by an entirely different method, also described in \cite[Ch. 3]{Vaughan}. One consequence of these results is that for any $p > 2$ one has
\begin{equation}\label{eq1} \Vert h_N \Vert_{L^p(\mathbb{T})} 
:= (\int_0^{1} |h_N(\theta)|^p\ d\theta)^{1/p} \ll_p N^{1 - 1/p} \log^{-1} N.\end{equation}
Disregarding the multiplicative constant, \eqref{eq1} is easily seen to be tight, this being a reflection of the prime number theorem, which in particular
gives $|h_N(\theta)| \sim N/\log N$ when $|\theta| \ll 1/N$.
In fact, one could estimate $\Vert h_N \Vert_p$ asymptotically.\vs

\ni In this paper we consider generalizations of this exponential sum in which the primes are replaced by such sets
as the twin primes.  More precisely, let $F = \prod_{j=1}^k (a_j n + b_j)$ be a polynomial which is product of $k \geq 1$ linear
factors with integer coefficients. We assume that $a_j \neq 0$ and that the discriminant
\[ \Delta := \prod_{1 \leq i < j \leq k} (a_i b_j - a_j b_i)\]
is non-zero, or in other words that no two linear factors are rational multiples of each other.
For any such $F$, we consider the set 
\[
X = X(F) := \{ n \in \Z^+: F(n) \hbox{ is the product of } k \hbox{ primes }\}.
\]
For example, if $F(n) := n(n+2)$, then $X$ is the collection of twin primes.
For any $q \geq 1$, we introduce the quantity $\gamma(q) = \gamma(q,F)$ defined by
\begin{equation}\label{rho-def}
\gamma(q) := q^{-1} |\{ n \in \Z/q\Z: (q, F(n)) = 1 \}|.
\end{equation} Observe that $0 \leq \gamma(q) \leq 1$. This number
represents the proportion of residue classes modulo $q$ which could hope to contain an infinite number of elements of $X$.
We will impose the natural \emph{non-degeneracy condition}
\begin{equation}\label{nondeg-def}
\gamma(q) > 0 \hbox{ for all } q \geq 1,
\end{equation}
since $X(F)$ is clearly finite if this condition fails.
\begin{remark} Note
that to verify \eqref{nondeg-def} it suffices by the Chinese remainder theorem to verify it for prime $q$, and in fact
one only needs to verify it for finitely many prime $q$, as the claim is trivial if $q$ is larger than both $k$
and the discriminant $\Delta$ of $F$.  Indeed for primes $p$, we have 
\begin{equation}\label{rho-explicit}
\gamma(p) = 1-k/p \hbox{ whenever } p > k \hbox{ and } p \nmid \Delta.
\end{equation}
\end{remark}
\ni The exponential sum over the twin primes, or more generally exponential sums of the form
\begin{equation}\label{eq2} 
h(\theta) = h_{N;F}(\theta) :=  \sum_{\substack{n \leq N \\ n \in X(F)}} e(n\theta),
\end{equation}
have not received much attention. One very good explanation for this is that, for a fixed choice of $F$, 
the set $X(F)$ is not even known to be non-empty. Thus there is certainly no hope of obtaining an asymptotic formula for $h$ on the major arcs, unless one is willing to obtain results on average over choices of $F$, as was done in \cite{Balog}. On the other hand, the \emph{Hardy-Littlewood $k$-tuple conjecture} asserts that for a fixed $F$
\begin{equation}\label{tuple}
h_{N;F}(0) = |X(F) \cap [1,N]| = (1 + o_{F}(1)) \mathfrak{S}_F \frac{N}{\log^{k} N}
\end{equation}
as $N \rightarrow \infty$,
where $\mathfrak{S}_F$ is the \textit{singular series}
\begin{equation}\label{gamma-def}
 \mathfrak{S}_F := \prod_p \frac{\gamma(p)}{(1 - \frac{1}{p})^k} \; .
\end{equation}
From \eqref{rho-explicit} we see that the infinite product in \eqref{gamma-def} is absolutely convergent. We will prove certain facts about $\mathfrak{S}_F$ later on, but to give the reader a feel for it we remark that if $|a_i|, |b_i| \leq N$ then $\mathfrak{S}_F \ll_{\epsilon} N^{\epsilon}$. Also, the $L^p$ average of $\mathfrak{S}_F$ over any sensible family of $F$ (for example $F(n) = n(n+h)$, $h = 1,\dots,N$) will be absolutely bounded.\vs

\ni The conjecture \eqref{tuple}
is still very far from being settled; however it is well known that sieve methods (either the Selberg sieve or the combinatorial sieve) give upper bounds of the form
\begin{equation}\label{eq3} |X(F) \cap [1,N]| \; \ll_k \; \mathfrak{S}_F(1 + o_{F}(1)) \frac{N}{\log^k N}\end{equation} for fixed $F$
(see e.g. \cite[Theorem 5.3]{Halberstam-Richert}). This upper bound only differs
from \eqref{tuple} by a multiplicative constant depending on $k$.  As an immediate consequence of this bound and 
Parseval's identity one may place an upper bound on $\Vert h_{N,F} \Vert_2$ which (conditionally on the $k$-tuple conjecture) is of the correct order of magnitude. Our first result (proven at the end of \S \ref{majorant-sec}) 
asserts that we in fact can obtain good estimates for all $L^p$ norms, $p>2$, of these exponential sums:

\begin{theorem}\label{mainthm}  Let $F$ be a product of $k$ integer linear forms having non-zero discriminant, obeying the non-degeneracy condition
\eqref{nondeg-def} and with coefficients $a_i,b_i$ satisfying $|a_i|,|b_i| \leq N$. Let $h_{N,F}$ be the exponential sum over $X(F) \cap [1,N]$ as 
defined in \eqref{eq2} above.  Then for every $2 < p < \infty$ we have the estimate
\[ \Vert h \Vert_p \; \ll_{p,k} \; \mathfrak{S}_F  N^{1 - 1/p}(\log N)^{-k}.\]
The implied constant depends on $k$ but is otherwise independent of $F$; it is explicitly computable.
\end{theorem}

\begin{remark}  In the converse direction, if the Hardy-Littlewood prime $k$-tuple conjecture
\eqref{tuple} is true, then one easily obtains the lower bound
\[
 \Vert h \Vert_p \; \gg \; \mathfrak{S}_F N^{1 - 1/p}(\log N)^{-k},
\]
at least for fixed $F$ and $N$ sufficiently large. This is a simple consequence of \eqref{tuple} and the observation that
$h_{N;F}(\theta) \geq \frac{1}{2} |X(F) \cap [1,N]|$ if $|\theta| \leq c_{k}/N$ and $c_{k} > 0$ is
sufficiently small.  However, when $p=2$, we see from Parseval's identity and the prime $k$-tuple conjecture
\eqref{tuple} that
$$ \Vert h \Vert_2 = |X(F) \cap [1,N]|^{1/2} \gg_k \mathfrak{S}_F^{1/2}  N^{1/2} (\log N)^{-k/2}$$
and so Theorem \ref{mainthm} is expected to fail logarithmically at the endpoint $p=2$. \end{remark}

\ni The proof of Theorem \ref{mainthm} proceeds not by studying the set $X(F) \cap [1,N]$ directly, but rather 
by constructing a slightly larger object $\beta_R$, which we call an \emph{enveloping sieve} for $X(F) \cap [1,N]$.  
This sieve function is considerably easier to analyze than $X(F) \cap [1,N]$ itself, and in particular
enjoys very good Fourier properties.  It then turns out that the only property we need on $X(F) \cap [1,N]$ itself to
deduce Theorem \ref{mainthm} is that it is \emph{majorized} by this sieve.  In particular, this theorem will also
hold if the set $X(F)$ is replaced by other similar sets, for instance if primes are replaced by almost primes, or
with any subset of the almost primes.  Combining this observation with an argument of the first 
author \cite{Green}, we prove the following result. 

\begin{theorem}\label{ap-theorem}
Define a \textit{Chen prime} to be a prime $p$ such that $p + 2$ is either prime or a product of two primes.
Then there are infinitely many $3$-term arithmetic progressions of Chen primes.
\end{theorem}

\begin{remark} Such a theorem, with $3$ replaced by a more general integer $k$, can probably be obtained by adapting the 
methods of \cite{green-tao-primes}, in which the primes were established to contain arithmetic progressions of any length $k$, though the explicit constants obtained by those methods are likely to be substantially worse.  
Our techniques here are somewhat similar in spirit to those in \cite{green-tao-primes}
in that they are based on establishing a \emph{transference principle}, in this case from Roth's theorem \cite{Roth}
concerning progressions of length three in a dense subset of the integers, to the corresponding result concerning
dense subsets of a sufficiently pseudorandom enveloping set.  In \cite{green-tao-primes} this transference principle 
was established (under quite strong pseudorandomness assumptions) 
by the finitary analogue of ergodic theory methods.  Here, we use a much simpler
harmonic analysis argument to prove the transference principle required for Theorem \ref{ap-theorem}, and as such
we require a much less stringent pseudorandomness condition (basically, we require the Fourier coefficients
of the enveloping sieve to be under control).  \end{remark}

\begin{remark} We note that results along the lines of Theorem \ref{ap-theorem} can be approached using sieve methods and the Hardy-Littlewood circle method in a more classical guise. For example, Tolev \cite{tolev} showed that there are infinitely many 3-term progressions $p_1 < p_2 < p_3$ of primes such that $p_i+2$ is a product of at most $r_i$ primes, where $(r_1,r_2,r_3)$ can be taken to be $(5,5,8)$ or $(4,5,11)$. \end{remark}

\ni The authors would like to thank the anonymous referee for a particularly attentive reading of the paper, and John Friedlander for a number of helpful remarks concerning the material in the appendix.

\section{Notation}\label{notation-sec}
\ni We take the opportunity to introduce here some basic notation common to the whole paper, most of which is standard in analytic number theory. Some of our notation, such as the use of the expectation symbol $\mathbb{E}$, is less standard in this field but is very convenient for certain additive problems of the type considered here. This is because there are ideas from ergodic theory bubbling just beneath the surface.\vs

\ni Let $k$ be a parameter. We write $A \ll_k B$ or $A = O_k(B)$ to denote the estimate $A \leq C_k B$ where $C_k > 0$ is a constant depending only on $k$. The notation $o_k(1)$ will be reserved for a quantity which tends to zero as the variable $N$ tends to $\infty$, at a rate which may depend on $k$ (and possibly other subscripted parameters or objects).\vs

\ni For a statement $P$, we will occasionally write $\mathbf{1}_P$ to denote the indicator of $P$, 
thus $\mathbf{1}_P = 1$ if $P$ is true and $\mathbf{1}_P = 0$ if $P$ is false.  If $A$ is a set, we use $\mathbf{1}_A$ to denote the function $\mathbf{1}_A(x) := \mathbf{1}_{x \in A}$. We will always write $|A| = \sum_x \mathbf{1}_A(x)$ for the cardinality of $A$.\vs

\ni We write $\Z_q$ for the cyclic group $\Z/q\Z$, and $\Z_q^* := \{ a \in \Z_q: (a,q) = 1 \}$ for the multiplicatively invertible elements of $\Z_q$. We will often identify $\Z_q$ with $\{0,\dots,q-1\}$ or with $\{1,\dots,q\}$ in the obvious manner.\vs

\ni As we remarked, it is convenient to work with the language of measures and conditional expectation in much the same way as the authors did in \cite{green-tao-primes}.
If $f: A \to \mathbb{C}$ is a function and $B$ is a non-empty finite subset of $A$,
we write $\E_{x \in B}f(x)$ for the average value of $f$ on $B$, that is to say
\[ \E_{x \in B}f(x) := \frac{1}{|B|} \sum_{x \in B} f(x).\] 
For the purposes of this paper this subscript notation is more convenient than the more traditional notation $\E(f(x) | x \in B)$ of probability theory.
In many cases we shall just have $B=A$, in which
case we abbreviate $\E_{x \in A}f(x)$ as $\E(f)$.\vs

\ni As remarked earlier, we write $e(\theta) = e^{2\pi i \theta}$. If $m$ is a positive integer then we introduce the related notation $e_m$ for the additive character on $\mathbb{Z}_m$ defined by $e_m(x) = e(x/m)$.

\section{The enveloping sieve}

\ni We now begin the proof of Theorem \ref{mainthm}.  The strategy of the sieve method is not to study the set $X \cap [1,N]$ directly, but rather to construct an \emph{enveloping sieve} $\beta_R: \Z^+ \to \R^+$ which is concentrated
on $X \cap [1,N]$ and for which one has satisfactory control of the Fourier coefficients of $\beta_R$.  We will normalize
$\beta$ to have average value comparable to 1; from \eqref{tuple}, we would thus expect $\beta_R$ to roughly be of size
$\log^{k} N / \mathfrak{S}_F$ on $X \cap [1,N]$, and this indeed will be the case; indeed, $\beta_R$
will be a \emph{pointwise majorant} of (most of) the normalized indicator function $c_{k}\mathfrak{S}_F^{-1} \log^k N  \mathbf{1}_{X \cap [1,N]}$.\vs

\ni Henceforth the polynomial $F$ will be fixed.  
We shall replace the
set $X = X(F)$ by ``localized'' versions $X_q$ for various moduli $q$, defined by
\[ X_q := \{ n \in \Z_q: (q, F(n)) = 1 \}.\]
This set has already appeared in \eqref{rho-def}.  Observe that when reduced $\md{q}$ all the elements in $X$, with at most $kq$
exceptions, lie in $X_q$, the exceptions arising from the cases when $a_jn+b_j \in \{1,\ldots,q\}$ for some
$1 \leq j \leq k$.\vs

\ni We can now describe the important properties of enveloping sieve, which is intended to majorize the set
\[ X_{R!} = \{ n \in \Z: (d,F(n)) = 1 \hbox{ for all } 1 \leq d \leq R \}\]
where $R$ is a large integer (the sieve parameter) to be chosen later. Eventually we will set 
$R := \lfloor N^{1/20} \rfloor$.  

\begin{proposition}[The enveloping sieve]\label{prop3} Let $F$ be the product of $k$ integer linear forms with non-zero discriminant,
obeying the non-degeneracy condition \eqref{nondeg-def} and with coefficients $a_i,b_i$ satisfying $|a_i|,|b_i| \leq N$. Let $R \leq N$ be a large integer.
Then there is a non-negative function $\beta := \beta_R: \mathbb{Z} \rightarrow \R^+$ 
with the following properties:
\begin{itemize}
\item[(i)] \textup{(Majorant property)} We have
\begin{equation}\label{sing-series} \beta(n) \gg_k \mathfrak{S}_F^{-1} \log^k R \mathbf{1}_{X_{R!}}(n)
\end{equation}
for all integers $n$. In particular, $\beta(n)$ is non-negative.  
\item[(ii)] \textup{(Crude upper bound)}  We have
\begin{equation}\label{upper-crude}
\beta(n) \ll_{k,\epsilon} N^{\epsilon}
\end{equation}
for all $0 < n \leq N$ and $\epsilon > 0$.
\item[(iii)] \textup{(Fourier expansion)} We have
\begin{equation}\label{eq3.20} \beta(n) \; = \; \sum_{q \leq R^2} \sum_{a \in \Z_q^*}w(a/q)e_q(-an),\end{equation}
where $w(a/q) = w_R(a/q)$ obeys the bound
\begin{equation}\label{star} |w(a/q)| \; \ll_{k,\epsilon} \; q^{\epsilon - 1}\end{equation}
for all $q \leq R^2$ and $a \in \Z^*_q$.  Also we have $w(0) = w(1) = 1$.
\item[(iv)] \textup{(Fourier vanishing properties)} Let $q \leq R^2$ and $a \in \Z_q^*$.  If $q$ is not square-free,
then $w(a/q) = 0$.  Similarly, if $\gamma(q) = 1$ and $q > 1$, then $w(a/q) = 0$.
\end{itemize}
\end{proposition}

\begin{remark} The enveloping sieve $\beta_R$ is essentially a normalised version of the one used by Selberg \cite{selberg} in obtaining the upper bound \eqref{eq3}. It will be constructed quite explicitly in \S \ref{appendix1}, where a full proof of Proposition \ref{prop3} will also be supplied. One does not need to know the exact construction in order to use Proposition
\ref{prop3} in applications. In \S \ref{appendix2} we offer some remarks which explain the phrase ``enveloping sieve'', as well as a comparison of $\beta_R$ (in the case $F(n) = n$) and a somewhat different majorant for the primes which we used in \cite{green-tao-primes}.\end{remark}

\section{The Hardy-Littlewood majorant property for the enveloping sieve}\label{majorant-sec}

\ni In Proposition \ref{prop3} we constructed a majorant $\beta_R$ for the set $X_{R!}$ which enjoyed good control
on the Fourier coefficients; indeed, from \eqref{star} it is easy to obtain good $L^p$ control on the Fourier 
coefficients of $\beta_R$  for $p>2$, although the bound in \eqref{star} causes unavoidable logarithmic losses at 
the endpoint $p=2$.  When $p=4,6,8,\ldots$ one can pass from this to good control on $h$ by Parseval's formula, 
and hence by interpolation one can obtain Theorem \ref{mainthm} in the case $p \geq 4$.  However this simple argument 
does not seem to suffice in the more interesting region\footnote{For the application to finding progressions of length three, we in fact need
this estimate for some $2 < p < 3$.} $2 < p < 4$ because of the logarithmic failure of the estimate at $p=2$ remarked on in the introduction.  
In this region $2 < p < 4$ (or more generally for $p \notin 2 \N$) there are no good monotonicity properties in $L^p$
of the Fourier transform to exploit; see \cite[p144]{Montgomery} for a simple example where monotonicity in this sense 
breaks down, and \cite{bachelis,green-ruzsa,mockenhaupt-schlag} for further discussion. \vs 

\ni On the other hand, all known examples where monotonicity breaks down are rather pathological. 
When the majorant $\beta_R$ enjoys additional Fourier or geometric properties it is often possible
to recover estimates of Hardy-Littlewood majorant type.  This is known as the \textit{restriction phenomenon}
and has been intensively studied in harmonic analysis; see for instance \cite{Tao} for a recent survey of this theory.\vs

\ni The application of ideas from restriction theory to number theory was initiated by Bourgain in the papers
\cite{Bourgain,bgafa} (among others). In \cite{Bourgain}, which is of the most relevance to us, it is shown how to obtain bounds for $\Vert \hat{f} \Vert_p$, where $f$ is a function supported on the primes. Such a bound was also obtained in a paper of the first author \cite{Green}, using the Brun sieve and somewhat precise versions of the prime number theorem for arithmetic progressions.  Our arguments here,
while certainly in a similar spirit, differ from the previous arguments in that we rely purely on the properties of
the majorant given by Proposition \ref{prop3}; in particular we do not need any type of control of the Riemann $\zeta$- function or of related objects.\vs

\ni Proposition \ref{rest-beta} gives a restriction estimate for $\beta$ in the setting of the finite abelian group $\mathbb{Z}_N$. The \textit{dual} form of this estimate, which is Proposition \ref{prop4}, will be of more interest for applications (indeed, Theorem \ref{mainthm} will be an almost immediate consequence of it). As in many arguments in restriction theory, the argument is based in spirit on the Tomas-Stein
method, first used in \cite{tomas}. 

\begin{proposition}[$L^q$ Restriction estimate for $\beta$]\label{rest-beta}  Let $R, N$ be large numbers such that $1 \ll R \ll N^{1/10}$ and let $k$, $F$, $\beta_R$ be as in Proposition \ref{prop3}.  Let $f: \Z_N \to \C$ be arbitrary.  Then for every $1 < q < 2$ we have
\begin{equation}\label{form1} \E_{1 \leq n \leq N}|\sum_{b \in \Z_N} f(b) e_N(bn)|^2 \beta_R(n)  \ll_{q,k} 
(\sum_{b \in \Z_N} |f(b)|^q)^{2/q}.\end{equation}
\end{proposition}

\begin{proof} We will first examine the left-hand side of \eqref{form1} in the special case when $f$ is supported on a set $B \subseteq \mathbb{Z}_N$ and $\Vert f \Vert_{\infty} \leq 1$. We assert that for every $\epsilon > 0$ we have
\begin{equation}\label{set-version} \E_{1 \leq n \leq N} |\sum_{b \in \mathbb{Z}_N} f(b) e_N(bn)|^2 \beta_R(n)  \ll_{\epsilon,k}  
|B|^{1+\epsilon}.\end{equation}
This statement, though at first sight weaker than \eqref{form1}, is actually equivalent to it. To deduce \eqref{form1} from \eqref{set-version}, normalise so that
\begin{equation}\label{eq444} (\sum_{b \in \Z_N} |f(b)|^q)^{2/q} = 1.\end{equation}
In particular this implies that $|f(b)| \leq 1$ for all $b \in \Z_N$.  If we then define the sets $B_j \subset \Z_N$
for all $j \geq 0$ as $B_j := \{ b \in \Z_N: 2^{-j-1} < |f(b)| \leq 2^{-j} \}$, then we have
$$ \sum_{b \in \Z_N} f(b) e_N(bn) = \sum_{j \geq 0} \sum_{b \in B_j} f(b) e_N(bn).$$
Also, from \eqref{eq444} we have $|B_j| \leq 2^{(j+1)q}$.  Applying \eqref{set-version} to
$B_j$ (and with $f(b)$ replaced by $2^j f(b)$) we have
\begin{eqnarray*} \big(\E_{1 \leq n \leq N} |\sum_{b \in B_j} f(b) e_N(bn)|^2 \beta_R(n)\big)^{1/2}  & \ll_{q,\epsilon,k} &
2^{-j} |B_j|^{(1+\epsilon)/2} \\ & \leq  & 2^{((1+\epsilon)q/2 - 1)j}\end{eqnarray*}
for any $\epsilon > 0$.  If we choose $\epsilon$ sufficiently small, then the exponent $(1+\epsilon)q/2 - 1$
is negative since $q < 2$.  Summing this in $j$ using the triangle inequality in $L^2_{\beta,N}$, defined as the space of complex sequences $(a_n)_{n =1}^N$ together with the inner product
\[ \langle (a_n), (b_n)\rangle = \E_{1 \leq n \leq N}a_n \overline{b_n}\beta_R(n),\] we see that \eqref{form1} indeed follows (since $q < 2$).\vs

\ni It remains to prove \eqref{set-version} which, recall, was to be demonstrated on the assumption that $\Vert f \Vert_{\infty} \leq 1$ and that $f$ is supported on $B$.
When $B$ is large, say $|B| \geq R^\epsilon$, we can use the crude estimate $\beta_R(n) \ll_{k,\epsilon} R^{\epsilon^2}$ (for instance) from \eqref{upper-crude}, combined with Parseval's identity 
to conclude the bound
\[ \E_{1 \leq n \leq N}|\sum_{b \in B} f(b) e_N(bn)|^2 \beta_R(n)  \ll_{k,\epsilon}  R^{\epsilon^2} |B| 
\leq |B|^{1+\epsilon}.\]
Thus we may take $|B| \leq R^\epsilon$.\vs

\ni We will essentially be doing Fourier analysis on $\mathbb{Z}_N$, but exponentials such as $e_q(a)$, which arise from the Fourier expansion \eqref{eq3.20}, will also appear. To facilitate this somewhat awkward juxtaposition
it will be convenient to introduce a smooth Fourier cutoff.  Let $\psi: \R \to \R^+$, $\Vert \psi \Vert_{\infty} \leq 1$, be a non-negative even bump function supported on the interval $[-1/10,1/10]$ whose Fourier transform $\widehat{\psi}(y) := \int_\R \psi(x) e(xy)\ dx$ is non-negative everywhere and bounded away from zero when $y \in [-1,1]$. Such a function can easily be constructed,
for instance by convolving a non-negative even bump function supported on $[-1/20, 1/20]$ with itself.  Then
we have
\[ \E_{1 \leq n \leq N} |\sum_{b \in B} f(b) e_N(bn)|^2 \beta_R(n)  \ll \frac{1}{N} \sum_{n \in \Z} |\sum_{b \in B} f(b)e_N(bn)|^2 \beta_R(n) \hat \psi(n/N).\]
We expand the right-hand side using the Fourier expansion of $\beta_R$, \eqref{eq3.20}, to obtain
\[ \frac{1}{N} \sum_{n \in \Z} \sum_{q \leq R^2} \sum_{a \in \Z_q^*} \sum_{b, b' \in B} f(b) \overline{f(b')}
w(a/q) e_N(bn) e_N(-b' n) e_q(-an) \widehat{\psi}(n/N).\]
Now the Poisson summation formula implies that
\[ \frac{1}{N} \sum_{n \in \Z} e_N(bn) e_N(-b'n) e_q(-an) \widehat{\psi}(n/N) = \sum_{m \in \Z} \psi\left(N(m + \frac{b-b'}{N} - \frac{a}{q})\right),\]
and hence by construction of $\psi$
\[ \left| \frac{1}{N} \sum_{n \in \Z} e_N(bn) e_N(-b'n) e_q(-an) \widehat{\psi}(n/N)\right| \leq \mathbf{1}_{\| \frac{b-b'}{N} - \frac{a}{q} \| \leq 1/N},\]
where $\|x\|$ denotes the distance of $x$ to the nearest integer.  Applying this and \eqref{star}, and
bounding the coefficients $f(b)$, $\overline{f(b')}$ by one, our task thus reduces to proving that
\begin{equation}\label{sum-to-bound} \sum_{q \leq R^2} \sum_{a \in \Z_q^*} \sum_{\substack{b,b' \in B\\ \| \frac{b-b'}{N} - \frac{a}{q} \| \leq 1/N }} q^{\epsilon/2-1} 
\ll_\epsilon |B|^{1+\epsilon}.\end{equation}
Let us first dispose of the contribution of the large $q$, for which $q \geq |B|$, where it is possible to estimate $q^{\epsilon/2-1}$ by $|B|^{\epsilon/2-1}$.  Observe that for any two distinct fractions $a/q$, $a'/q'$ with $q, q' \leq R^2$
and $a \in \Z_q^*$, $a' \in \Z_{q'}^*$, we have
\[ \left\|\frac{a}{q} -\frac{a'}{q'}\right\| \geq \frac{1}{qq'} \geq \frac{1}{R^4} \geq \frac{2}{N}\]
since $N \gg R^{10}$.  This shows that for each fixed $b, b'$ there is at most one fraction $\frac{a}{q}$ which contributes to the sum \eqref{sum-to-bound}.  Thus the contribution of the large $q$ is at most $|B|^2 |B|^{\epsilon/2-1}$, which is acceptable.  It
will thus suffice to estimate the sum over small $q$, that is to say those $q$ for which $q \leq |B|$, and to that end it is enough to prove that 
\begin{equation}\label{sum-to-bound2} \sum_{q \leq |B|} \sum_{a \in \Z_q^*} \sum_{\substack{b,b' \in B\\ \| \frac{b-b'}{N} - \frac{a}{q} \| \leq 1/N }} q^{-1} 
\ll_\epsilon |B|^{1+\epsilon/2}.\end{equation}
We estimate the left-hand side by
\[  c\sum_{q \leq |B|} \sum_{a \in \Z_q} \sum_{b,b' \in B} q^{-1} 
\sum_{m \in \Z} \widehat{\psi}\left(N(m - \frac{b-b'}{N} + \frac{a}{q})\right),\]
where $c = \max_{y \in [-1,1]} |\widehat{\psi}(y)|^{-1}$,
which by the Poisson summation formula is
\[ cN^{-1}\sum_{q \leq |B|} \sum_{a \in \Z_q} \sum_{b,b' \in B} q^{-1} 
\sum_{n \in \Z}  \psi(n/N) e_N((b-b')n) e_q(-an).\]
Performing the summation over $a$, this becomes
\[  cN^{-1}\sum_{q \leq |B|} \sum_{b,b' \in B} 
\sum_{\substack{n \in \Z \\ q | n}} \psi(n/N) e_N((b-b')n)\]
which we can rearrange as
\begin{equation}\label{sum-to-bound3}  cN^{-1} \sum_{n \in \Z} \psi(n/N) \left|\sum_{b \in B} e_N(bn)\right|^2d_{|B|}(n),\end{equation}
where
\[ d_{|B|}(n) := \sum_{\substack{q \leq |B| \\ q | n}} 1.\]
The $n=0$ term of this sum may be bounded, very crudely, by $O(N^{-1}|B|^3 )$. Since $|B| \leq R^{\epsilon} \ll N^{1/10}$, this term is certainly at most $O(|B|)$ and hence does not make an important contribution to \eqref{sum-to-bound2}.
Let us thus discard the $n=0$ term, and bound the remaining terms using H\"older's inequality.  Since $\psi(n/N)$
is supported in the region $|n| \leq N/2$, we thus obtain
\[ \eqref{sum-to-bound3} \ll \big(\E_{0 < |n| \leq N/2} |\sum_{b \in B} e_N(bn)|^{2+\epsilon}\big)^{\frac{2}{2+\epsilon}}
\big(\E_{0 < |n| \leq N/2} d_{|B|}(n)^{\frac{2 + \epsilon}{\epsilon}}  \big)^{\frac{\epsilon}{2+\epsilon}}.\]
Parseval and the crude bound $|\sum_{b \in B} e_N(bn)| \leq |B|$ yield
\[ \big(\E_{0 < |n| \leq N/2} |\sum_{b \in B} e_N(bn)|^{2+\epsilon} \big)^{\frac{2}{2+\epsilon}}
\leq (|B|^{1+\epsilon})^{\frac{2}{2+\epsilon}} = |B|^{1 + \frac{\epsilon}{2+\epsilon}}.\]
On the other hand, from standard moment estimates for the restricted divisor function 
$d_{|B|}$ (see for instance \cite{Bourgain,Ruzsa}) we have
\[ \E_{0 < |n| \leq N/2}d_{|B|}(n)^m  \ll_{\epsilon,m} |B|^\epsilon\]
for any $m \geq 1$ (in fact one can replace $|B|^\epsilon$ with $\log^{2^m-1} |B|$ when $m$ is an integer - see \cite{Ruzsa}).
Applying this with $m := \frac{2 + \epsilon}{\epsilon}$ and working back, one gets $\eqref{sum-to-bound3} \ll_{\epsilon} |B|^{1 + \epsilon/2}$, which in turn implies \eqref{sum-to-bound2} and \eqref{sum-to-bound}, and hence the proposition.
\end{proof}\vs

\ni As promised, we now give what is essentially a dual form of Proposition \ref{rest-beta}.

\begin{proposition}[$L^p$ extension estimate for $\beta$]\label{prop4} Let $R, N$ be large numbers such that $1 \ll R \ll N^{1/10}$ and let $k$, $F$, $\beta_R$ be as in Proposition \ref{prop3}. Let $\{a_n\}_{n \leq N}$ be an arbitrary sequence of complex numbers, and let $p > 2$ be a fixed real number. Then we have
\begin{equation}\label{fixed-rot}
\big(\sum_{b \in \Z_N} \left| \E_{1 \leq n \leq N}a_n \beta_R(n) e_N(-bn)  \right|^p\big)^{1/p} \ll_{p,k}  \big(\E_{1 \leq n \leq N}  |a_n|^2 \beta_R(n)\big)^{1/2}
\end{equation}
and
\begin{equation}\label{variable-rot}
 \left\Vert \E_{1 \leq n \leq N} a_n \beta_R(n) e(n\theta)  \right\Vert_{L^p(\mathbb{T})} \; \ll_{p,k} \;  N^{-1/p} \big(\E_{1 \leq n \leq N} |a_n|^2 \beta_R(n) \big)^{1/2}.
\end{equation}
\end{proposition}

\begin{proof}
The first claim follows from Proposition \ref{rest-beta} and a duality argument. We will expand on this point somewhat using the language of operators, since this will afford us an opportunity to at least partially explain the terms \textit{restriction} and \textit{extension}. Everything we say could also be phrased in terms of inequalities on exponential sums. Recall that $L^2_{\beta,N}$ is the space of complex sequences $(a_n)_{n =1}^N$ together with the inner product
\[ \langle (a_n), (b_n)\rangle = \E_{1 \leq n \leq N}a_n \overline{b_n}\beta_R(n) .\] Consider also the space of complex-valued functions on $\mathbb{Z}_N$ together with the $l^q$ norms
\[ \Vert f \Vert_{l^q(\Z_N)} := (\sum_{x \in \Z_N} |f(x)|^q)^{1/q}\] and the inner product $\langle f, g \rangle := \sum_{x \in \Z_N} f(x) \overline{g(x)}$. Here, $\{1,\dots,N\}$ may perhaps be thought of as a finite version of $\mathbb{Z}$, whilst $\mathbb{Z}_N$ plays the r\^ole of a discretised circle $\mathbb{T}$. Now define the \textit{restriction} map $T : l^q(\Z_N) \rightarrow L^2_{\beta,N}$ by \[ f \mapsto (\sum_{x \in \Z_N} f(x)e_N(nx))_{n=1}^N\] (the reason for the name is that $T$ can be thought of as the restriction of the discrete Fourier transform on $\mathbb{Z}_N$ to the weight $\beta$). One can check that the adjoint of $T$ is the \textit{extension} map $T^* : L^2_{\beta,N} \rightarrow l^p(\Z_N)$ defined by \[(a_n)_{n=1}^N \mapsto \E_{1 \leq n \leq N}a_n \beta_R(n) e_N(-nx) .\] Here, of course, $p$ is the dual exponent to $q$ defined by the relation $p^{-1} + q^{-1} = 1$; note, however, that all of the spaces $l^r(\mathbb{Z}_N)$ are equivalent as vector spaces.\vs

\ni Now one can check that \eqref{form1}, the main result of Proposition \ref{rest-beta}, is equivalent to \[\Vert T \Vert_{l^q(\Z_N) \to L^2_{\beta,N}} \ll_{p,k} 1.\] On the other hand \eqref{fixed-rot} is equivalent to the statement 
\[ \Vert T^* \Vert_{L^2_{\beta,N} \to l^p(\Z_N)} \ll_{p,k} 1.\]
By the duality principle for operator norms one sees that these two statements are completely equivalent, and so \eqref{fixed-rot} is an immediate deduction from Proposition \ref{rest-beta}.\vs

\ni To prove \eqref{variable-rot}, observe that for any $0 \leq \theta < 1$
we can obtain the estimate
\[
 \sum_{b \in \Z_N} \left| \E_{1 \leq n \leq N} a_n \beta_R(n) e((b+\theta)n/N) \right|^p \; \ll_{p,k} \; 
\big(\E_{1 \leq n \leq N}|a_n|^2 \beta_R(n) \big)^{p/2} \]
from \eqref{fixed-rot} simply by multiplying each coefficient $a_n$ by $e_N(\theta n)$.  Integrating this in $\theta$ from 0 to 1 gives \eqref{variable-rot}.\end{proof} 

\begin{remark} A result of Marcinkiewicz \cite{marcinkiewicz-zygmund} allows one to obtain \eqref{fixed-rot} from \eqref{variable-rot}; see \cite{Green} for further discussion.\end{remark}

\ni The utility of Proposition \ref{prop4} is markedly increased by the next lemma.
\begin{lemma}\label{beta-L1}
Suppose that $R \leq \sqrt{N}$. Then
\[ \E_{1 \leq n \leq N}\beta_R(n)  \ll 1.\]
\end{lemma}
\begin{proof}
We make further use of the smooth cutoff function $\psi$. We have \[  \E_{1 \leq n \leq N} \beta_R(n)  \ll \frac{1}{N} \sum_{n \in \Z} \hat{\psi}(n/N) \beta_R(n).\]
Substituting in the Fourier expansion \eqref{eq3.20} of $\beta_R$ and using the Poisson summation formula, we have
\[ N^{-1}\sum_{n \in \Z} \hat \psi(n/N) \beta_R(n) = 
\sum_{q \leq R^2} \sum_{a \in \Z_q^*} w(a/q) \sum_{m \in \Z} \hat \psi\big(N(m-\frac{a}{q})\big).\]
Observe that since $R^2 \leq N$, the only time that $\hat \psi(N(m - \frac{a}{q}))$ is non-zero is when $m=a=0$ and $q=1$.
Thus, by \eqref{star}, we obtain
\[ \frac{1}{N} \sum_{n \in \Z} \widehat{\psi}(n/N) \beta_R(n) = w(0) \psi(0) \ll 1,\]
and the claim follows. \end{proof}\vs

\ni Observe that the last lemma, combined with Proposition \ref{prop3}(i), implies the classical sieve estimate \eqref{eq3}. With these preliminaries in place we can now quickly conclude the proof of Theorem \ref{mainthm}.\vs

\ni\textit{Proof of Theorem \ref{mainthm}.}  Fix $R := \lfloor N^{1/20} \rfloor$ (for instance).  Without loss of generality
we may assume $N$ to be sufficiently large depending on $k$ and $p$.  We now apply Proposition \ref{prop4}
to the sequence $a_n := \mathbf{1}_{X \cap X_{R!}}(n)/\beta(n)$.  From \eqref{sing-series}, \eqref{variable-rot} and Lemma \ref{beta-L1} we easily obtain 
\[ \big\Vert \sum_{X \cap X_{R!} \cap [1,N]} e(n\theta) \big\Vert_{L^p(\mathbb{T})} \; \ll_{p,k} \; 
\frac{\mathfrak{S}_F}{\log^k R} N^{1 - 1/p}.\]
Since all but at most $kR$ elements of $X$ lie in $X_{R!}$, the claim now follows since $R \sim N^{1/20}$ and $N$
is assumed sufficiently large.
\endproof

\begin{remark}  A straightforward modification of the above argument also gives the estimate
\[ \int_0^{1}
|\sum_{n \leq N} (\prod_{j=1}^k \Lambda(a_jn+b_j)) e(n\theta) |^p\ d\theta \; \ll_{p,k} \; 
\mathfrak{S}_F^p N^{p - 1},\]
where $\Lambda$ is the von Mangoldt function.\end{remark}

\section{A transference principle for 3-term arithmetic progressions}

\ni There is a well-known theorem of Roth \cite{Roth}, partially answering a question of Erd\H{o}s and Tur\'an, which states that any subset of $\{1,\dots,N\}$ with size at least $cN/\log\log N$ contains an arithmetic progression of length 3. In this section we will prove a result in this spirit for sets which are majorised by a ``pseudorandom measure''. In \S \ref{sec7} we will apply this result to deduce that there are infinitely many 3-term progressions of Chen primes (primes $p$ such that $p + 2$ is prime or a product of two primes).\vs

\ni We refer to the result (Proposition \ref{harmonic-transference-prop}) as a transference principle because, rather than prove it from first principles, we deduce it from Roth's theorem. For the reader familiar with \cite{green-tao-primes}, we remark that the transference principle we prove here is shown to hold under harmonic analysis conditions which should be regarded as much less stringent than the combinatorial conditions of that paper. However, our arguments here only apply to progressions of length 3. Indeed, it is by now well-understood that traditional harmonic analysis arguments cannot suffice to deal with progressions of length 4 or more (see \cite{gowers-4} for a further discussion).\vs

\ni We begin by recalling a formulation of Roth's theorem \cite{Roth} due to Varnavides \cite{varnavides}.

\begin{theorem}\label{bourg-th} \cite{Roth,varnavides}  Let $0 < \delta \leq 1$ be a positive numbers, and
let $N$ be a large prime parameter. Suppose that $f : \Z_N \rightarrow \mathbb{R}^+$ is a function satisfying the 
uniform bounds
\[ 0 \leq f(n) \leq 1 \hbox{ for all } n \in \Z_N\]
and the mean property
\[ \E(f) \geq \delta.\]
Then we have
\begin{equation}\label{eq6.1} \E_{n,d \in \Z_N}f(n)f(n+d)f(n+2d) \geq c(\delta)\end{equation}
for some constant $c(\delta) > 0$ depending only on $\delta$.
\end{theorem}

\begin{remark} By combining Bourgain's refinement \cite{bourg} of Roth's theorem with the argument of Varnavides \cite{varnavides} one can obtain the more explicit value
$c(\delta) = \exp(-C \delta^{-2} \log(1 + \frac{1}{\delta}))$ for $c(\delta)$; see also \cite{Green} for further discussion.  However, our arguments will not be sensitive to the exact value of $c(\delta)$.\end{remark}

\ni We now generalise the above proposition to the case where $f$ is not bounded by 1, but is instead bounded by some
suitably pseudorandom function.

\begin{definition}
We shall normalize the Fourier transform $\widehat{f}: \Z_N \to \C$ of a function $f: \Z_N \to \C$ according to the formula
\[ \widehat{f}(a) := \E_{n \in \Z_N}f(n) e_N(an)  \hbox{ for all } a \in \Z_N,\]
 and recall the \emph{Fourier inversion formula}
\begin{equation}\label{Fourier-inversion}
f(n) = \sum_{a \in \Z_N} \widehat{f}(a) e_N(-an).
\end{equation}
With this normalisation it is natural to measure the Fourier transform $\widehat{f}$ using the $l^p$ norms
\[ \Vert \widehat{f} \Vert_{p}^p := \sum_{a \in \Z_N} |\widehat{f}(a)|^p\]
for $1 < p < \infty$, with the usual convention that $\| \widehat{f} \|_{\infty} := \sup_{a \in \Z_N} |\widehat{f}(a)|$.
\end{definition}

\begin{proposition}\label{harmonic-transference-prop}
Let $N$ be a large prime parameter.
Suppose that $\nu : \Z_N \rightarrow \mathbb{R}^+$ and $f : \Z_N \rightarrow \mathbb{R}^+$ are non-negative
functions satisfying the majorization condition
\begin{equation}\label{major}
 0 \leq f(n) \leq \nu(n) \hbox{ for all } n \in \Z_N 
\end{equation}
and the mean condition
\begin{equation}\label{mean}
 \mathbb{E}_{n \in \Z_N}f(n) \geq \delta
\end{equation}
for some $0 < \delta \leq 1$. Suppose that $\nu$ and $f$ also obey the pseudorandomness conditions
\begin{equation}\label{nu-random}
|\widehat \nu(a) - \delta_{a,0}| \leq \eta \hbox{ for all } a \in \Z_N 
\end{equation}
where $\delta_{a,0}$ is the Kronecker delta, and
\begin{equation}\label{f-random}
\Vert \widehat{f} \Vert_{q} \leq M
\end{equation}
for some $\eta, M > 0$ and some $2 < q < 3$.  Then we have
\begin{equation}\label{f-conclude}
 \E_{n,d \in \Z_N}f(n)f(n+d)f(n+2d)  \geq \textstyle\frac{1}{2}\displaystyle c(\delta) 
- O_{\delta,M,q}(\eta)
\end{equation}
where $c(\delta) > 0$ is the same constant that appears in Theorem \ref{bourg-th}.
\end{proposition}

\begin{remark} Observe that the $\nu \equiv 1$, $\eta = 0$ case of this proposition is basically Theorem \ref{bourg-th}, with
an (inconsequential) loss of $\frac{1}{2}$ on the right-hand side of \eqref{eq6.1}.  This theorem should be compared 
with \cite[Theorem 3.5]{green-tao-primes}.  In practice, the condition \eqref{f-random} seems to be most easily verified
by first establishing a Hardy-Littlewood majorant property such as Proposition \ref{prop4} for \emph{all} functions $f$
majorized by $\nu$, not just the specific function $f$ of interest.\end{remark}

\begin{proof} We repeat the arguments of \cite[\S 6]{Green}.  Let
$0 < \epsilon \ll 1$ be a small parameter to be chosen later, and let
\begin{equation}\label{omega-def}
\Omega := \{a \in \Z_N : |\widehat{f}(a)| \geq \epsilon \}
\end{equation}
be the set of large Fourier coefficents of $f$.  
Write $B := B(\Omega,\epsilon)$ for the \textit{Bohr set}
\begin{equation}\label{bohr-def}
 B := B(\Omega,\epsilon) := \{ m \in \Z_N : |1 - e_N(am)| \leq \epsilon \;\; \mbox{for all } a \in \Omega\}.
\end{equation}
We shall decompose $f = f_1 + f_2$, where $f_1$ is the \emph{anti-uniform} component of $f$, defined by
\begin{equation}\label{f1-def} f_1(n) := \E_{m_1,m_2 \in B} f(n + m_1 - m_2) ,
\end{equation}
and $f_2 := f - f_1$ will be the \emph{uniform} component.  We can thus split the left-hand side of \eqref{f-conclude} into eight components
\begin{equation}\label{component} \E_{n,d \in \Z_N}f_i(n)f_j(n+d)f_k(n+2d) 
\end{equation}
where $i,j,k \in \{1,2\}$.\vs

\ni Let us first consider the term \eqref{component} corresponding to $(i,j,k) = (1,1,1)$, which will be the
dominant term.  The pseudorandomness condition \eqref{nu-random} can be used to obtain 
uniform bounds on $f_1$. Indeed using that together with \eqref{Fourier-inversion}, \eqref{major} and Parseval's
identity we obtain
\begin{align*}
|f_1(n)| &= |\E_{m_1,m_2 \in B}f(n + m_1 - m_2)|  \\ &\leq \E_{m_1,m_2 \in B} \nu(n + m_1 - m_2) \\ 
&= |\E_{m_1,m_2 \in B} \sum_{a \in \Z_N} \widehat{\nu}(a) e_N(-an-am_1+am_2) | \\
&= |\sum_{a \in \Z_N} \hat \nu(a) e_N(-an) |\E_{m \in B} e_N(-am) |^2| \\
&\leq \sum_{a \in \Z_N} |\hat \nu(a)| |\E_{m \in B} e_N(-am)|^2  \\
&\leq 1 + \eta \sum_{a \in \Z_N} |\E_{m \in B} e_N(-am) |^2 \\
&= 1 + \eta N / |B|.
\end{align*}
Note, however, that by the pigeonhole principle we have
\[ |B| \geq (c\epsilon)^{|\Omega|} N\]
for some absolute constant $c > 0$, while from \eqref{f-random} and Chebyshev's inequality we have
\[
|\Omega| \leq  (M/\epsilon)^q.
\]
Finally, $f_1$ is manifestly non-negative by \eqref{f1-def} and \eqref{major}.
We thus conclude that
\[
0 \leq f_1(n) \leq 1 + O_{\epsilon, M,q}(\eta) \hbox{ for all } n \in \Z_N.
\]
Furthermore from \eqref{mean} we have 
\[\E(f_1) = \E(f) \geq \delta.\]
Applying Theorem \ref{bourg-th} (possibly modifying $f_1$ by $O_{\epsilon,M,q}(\eta)$ first)
we obtain
\begin{equation}\label{f1-aps} 
\E_{x,d \in \Z_N}f_1(x)f_1(x+d)f_1(x+2d)  \geq 
c( \delta) - O_{\epsilon, \delta, M,q}(\eta).
\end{equation}
Now we consider the other seven terms of the form \eqref{component}, when at least one of the $i,j,k$ is equal to 2. The treatment of the seven cases is essentially the same, so for simplicity we consider only $i = j = k = 2$. Observe that for all $a \in \Z_N$ we have
\begin{equation}\label{fft-form}
\begin{split}
\widehat{f_1}(a) &= \E_{n \in \Z_N} \E_{m_1,m_2 \in B} f(n + m_1 - m_2)  e_N(an)  \\
&= \E_{n \in \Z_N} \E_{m_1,m_2 \in B} f(n+m_1-m_2) e_N(a(n+m_1-m_2) \times \\ &  \qquad \times \; e_N(-am_1) e_N(am_2)  \\
&= \widehat{f}(a) |\E_{m \in B} e(-am) |^2,
\end{split}
\end{equation}
and thus
\begin{equation}\label{eq445} | \widehat{f_2}(a)|  =  |\widehat{f}(a)| (1 - |\E_{m \in B}e_N(-am) |^2) \leq |\widehat{f}(a)|.\end{equation} This can be analysed in two different ways according as $a$ does or does not lie in $\Omega$.
When $a \not \in \Omega$ we have
\[ | \widehat{f_2}(a)| \leq |\widehat{f}(a)| \leq \epsilon\]
by \eqref{omega-def}.  On the other hand, when $a \in \Omega$, then $|e_N(-am) - 1| \leq \epsilon$ for all $m \in B$ by definition \eqref{bohr-def},
and hence in this case we have
\[ | \widehat{f_2}(a)| = |\widehat{f}(a)| (1 - |\E_{m \in B} e_N(-am) |^2) \leq 
3\epsilon\E(\nu) \leq 3\epsilon(1 + \eta)\]
thanks to the $a=0$ case of \eqref{nu-random}. Thus certainly
\begin{equation}\label{f2-Linfty} \| \widehat{f_2} \|_{\infty} = O_{\eta}(\epsilon).\end{equation}
Now by orthogonality we obtain the well-known identity
\[
\E_{n,d \in \Z_N}f_2(n)f_2(n+d)f_2(n+2d)  = \sum_{a \in \Z_N} \widehat{f_2}(a)^2 \widehat{f_2}(-2a).
\]
Such an expression is often treated using an $L^2$--$L^{\infty}$ inequality, but for us this is too expensive since we lack good $L^2$ control on $\nu$ (and hence on $f_2$).  Instead,
we shall use the $L^{q}$ hypothesis \eqref{f-random}.  By H\"older's inequality we obtain
\begin{equation}\label{eq200} |\E_{n,d \in \Z_N} f_2(n)f_2(n+d)f_2(n+2d) | \leq 
\| \widehat{f_2} \|_{q}^{q} \| \widehat{f_2} \|_{\infty}^{3-q}.\end{equation}
From \eqref{f-random} and \eqref{eq445} we have
\[\| \widehat{f_2} \|_{q} \leq \| \widehat{f} \|_q \leq M.\] Recalling \eqref{f2-Linfty} we therefore obtain, by \eqref{eq200}, the bound
\[ \E_{n,d \in \Z_N}f_2(n)f_2(n+d)f_2(n+2d) = O( M^q \epsilon^{3-q} ).\]
By choosing $\epsilon$ sufficiently small depending on $\delta$, $M$, $q$, this can be made less than $\frac{1}{14} c(\delta)$.  The same goes for the six other terms of the form \eqref{component} with $(i,j,k) \neq (1,1,1)$. Combining this with \eqref{f1-aps}, the claim follows.
\end{proof}

\section{Arithmetic progressions of Chen primes} \label{sec7}

\ni In this section we prove Theorem \ref{ap-theorem}, which stated that there are infinitely many 3-term progressions of Chen primes. A Chen prime, recall, is a prime $p$ such that $p + 2$ is a product of at most two primes, and Chen's famous theorem \cite{Chen} is that there are infinitely many such primes. Chen actually proved a somewhat stronger result, in which the smallest prime factor of $p + 2$ can be bounded below by $p^{1/10}$. It is important for us to have this extra information (though any positive number would do in place of $1/10$). In Iwaniec's unpublished notes \cite{iwaniec-notes} one may find a proof of Chen's theorem which leads to the value $\gamma = 3/11$. For the purposes of this section we will use the term Chen prime to refer to a $p$ for which $p + 2$ is either prime or a product $p_1p_2$, with $p_1,p_2 > p^{3/11}$. We quote, then, the following result from \cite{iwaniec-notes}.\footnote{Theorem \ref{chen-thm} is stronger than the result stated in \cite{iwaniec-notes} in one tiny way - Iwaniec shows that there are many Chen primes in $[0,N]$, not $(N/2,N]$. However it is completely clear that the proof can be adapted in a simple way to cover this case. Alternatively, a pigeonhole argument could be used to show that $(N/2,N]$ contains lots of Chen primes for many $N$, which would suffice for our purposes. } 

\begin{theorem}\label{chen-thm}  Let $N$ be a large integer.  Then the 
number of Chen primes in the interval $(N/2,N]$ is at least $c_1 N/\log^2 N$, for some absolute 
constant $c_1 > 0$. 
\end{theorem}

\ni For our argument we do not need to know an exact value for $c_1$.  
To have a chance of deducing Theorem \ref{ap-theorem} from Proposition \ref{harmonic-transference-prop}, we shall 
need a suitably pseudorandom measure $\nu$ on the cyclic group $\Z_N$ which majorises the Chen primes on $[N/4,N/2]$ in a reasonably efficient way. Unfortunately, the Chen primes are so irregularly distributed 
in progressions $a \pmod{q}$ ($q$ small) that there is no $\nu$ which would obey the pseudorandomness condition
\eqref{nu-random} for an $\eta$ which is small enough to be useful.  To get around this difficulty, we shall employ a device which we call the $W$-trick.  Let $t\gg 1$ be a very large real number (independent of $N$) to be specified later, 
and write $W = W_t := \prod_{p \leq t} p$.   Let $X_W \subseteq \Z_W$ denote the residue classes $b \in \mathbb{Z}_W$ such that $b$ and $b+2$ are both coprime to $W$; observe that
\begin{equation}\label{s-comp}
 |X_W| = W \prod_{3 \leq p \leq t} (1 - \frac{2}{p}) \ll |W| / \log^2 t.
\end{equation}
Note that all but $O(W)$ of the Chen primes lie in one of the residue classes associated to $X_W$.  From Theorem \ref{chen-thm} we thus have
$$ \sum_{b \in X_W} |\{ N/4 \leq n \leq N/2 : Wn + b \hbox{ is a Chen prime}\}| \gg \frac{WN}{\log^2(WN)} - O(W).$$
If we assume $N$ sufficiently large depending on $W$ (and hence on $t$), the right-hand side is $\gg WN/\log^2 N$.
From \eqref{s-comp} and the pigeonhole principle we may therefore choose $b \in X_W$ such that
\begin{equation}\label{lots-of-chens} 
|X| \gg \frac{N \log^2 t}{\log^2 N},
\end{equation}
where $X$ is the set
$$ X := \{ N/4 \leq n \leq N/2 : Wn + b \hbox{ is a Chen prime}\}.$$
\ni Fix $b \in X_W$ as above, and consider the polynomial $F(n) := (Wn + b)(Wn + b + 2)$.
To compute $\mathfrak{S}_F$, we observe that $\gamma(p) = 1$ when $p \leq t$ and $\gamma(p) = 1-2/p$ otherwise.  Since
the Euler product $\prod_p (1-2/p)/(1-1/p)^2$ is convergent, we see by the same computation used in 
\eqref{s-comp} that
\begin{equation}\label{sing-series-comp}
\log^2 t \ll \mathfrak{S}_F \ll \log^2 t.
\end{equation}
\ni Now let $R := \lfloor N^{1/20} \rfloor $ and let $\beta_R : \Z \rightarrow \mathbb{R}^+$ be the enveloping sieve associated to $F$ by Proposition \ref{prop3}.  We restrict this function to the set $\{1,\dots,N\}$, which we identify with $\Z_{N}$, and let $\nu: \Z_{N} \to \R^+$ be the resulting function.

\begin{lemma}\label{w-trick-lem}  For all $a \in\Z_{N}$ and $\epsilon > 0$, we have
\[ \hat \nu(a) = \delta_{a,0} + O_{\epsilon} (t^{\epsilon-1})\]
where $\delta_{a,0}$ is the Kronecker delta function.
\end{lemma}

\begin{proof} From \eqref{eq3.20} we have
\[
\widehat{\nu}(a) = \E_{n \in \Z_N} \beta_R(n) e_N(an)  = \sum_{q \leq R^2} \sum_{\substack{b \in \Z_q\\ (b,q) = 1}}
w(b/q)\E_{n \in \Z_N} e_N(an) e_q(-bn) .\]
When $1 < q \leq t$ we have $\gamma(q) = 1$, and hence $w(b/q) = 0$ by Proposition \ref{prop3}(iv).  By \eqref{star} we thus have
\begin{equation}\label{v-est}
\widehat{\nu}(a) = \E_{n \in \Z_N} e_N(an) +
O_\epsilon(t^{\epsilon-1} \sum_{t \leq q \leq R^2} \sum_{\substack{b \in \Z_q\\ (b,q) = 1}}
|\E_{n \in \Z_N} e_N(an) e_q(-bn) |).\end{equation}

\ni Observe that $\E_{n \in \Z_N} e_N(an)$ is just $\delta_{a,0}$.  From standard estimates
on the Dirichlet kernel we have
\begin{equation}\label{dirich-est} |\E_{n \in \Z_N} e_N(an) e_q(-bn)| \ll \min( 1, \frac{1}{N \| \frac{a}{N} - \frac{b}{q} \|} ).\end{equation}
All the fractions $b/q$, $q \leq R^2$, are separated from each other by at least 
$1/R^4 \leq 1/N$.  Thus we have $1/N \| \frac{a}{N} - \frac{b}{q} \| 
\ll R^4/N$ with at most one exception, which leads in view of \eqref{v-est} and \eqref{dirich-est} to \[
 \hat\nu(a) = \delta_{a,0} + O_\epsilon(t^{\epsilon-1} ( 1 + R^4 \frac{R^4}{N} )).\]
The claim follows.
\end{proof}\vs

\noindent\textit{Proof of Theorem \ref{ap-theorem}.} We fix $t$ to be an absolute constant to be chosen later, 
and take $N$ to be a sufficiently large prime
depending on $t$.
We define the function $f: \Z_N \to \R^+$ by
\[ f(n) := c \frac{\log^2 N}{\log^2 t} \mathbf{1}_X(n),\] and verify the various conditions of Proposition \ref{harmonic-transference-prop} for $f$.\vs

\ni First of all observe from Proposition \ref{prop3}(i) and \eqref{sing-series-comp} that we have the majorization
property \eqref{major} if $c > 0$ is chosen sufficiently small. Secondly, from \eqref{lots-of-chens} we have the condition \eqref{mean} with some absolute constant $\delta > 0$ independent of $N$ or $t$, as long as $N$ is large enough.  From Proposition \ref{prop4} (specifically \eqref{fixed-rot}) and Lemma \ref{beta-L1} we obtain the bound $\Vert \widehat{f} \Vert_{5/2} \leq M$ for some $M > 0$ independent of $N$ or $t$, if $N$ is sufficiently large depending on $t$. This confirms condition \eqref{f-random} (with $q = 5/2$).\vs

\ni Finally, from Lemma \ref{w-trick-lem} we obtain \eqref{nu-random} 
with $\eta = O_\epsilon(t^{\epsilon-1})$ and $M_0 = 1 + O_\epsilon(t^{\epsilon-1})$. \vs

\ni Thus we may indeed apply Proposition \ref{harmonic-transference-prop}, and doing so we obtain
\[ \E_{x,d \in \Z_N}f(x)f(x+d)f(x+2d)  \geq c - O_\epsilon(t^{\epsilon-1})\] for some 
absolute constant $c > 0$. By choosing $\epsilon := 1/2$ (say) and $t$ sufficiently large, we thus obtain
\[ \E_{x,d \in \Z_N}f(x)f(x+d)f(x+2d)  \geq c/2.\]
From the definition of $f$, this yields
$$ |\{ (x,d) \in \Z_N^2 : x, x+d, x+2d \in X \}| \gg_t N^2 / \log^6 N.$$
The degenerate case $d=0$ yields at most $O(N)$ pairs and can be discarded when $N$ is large.  Since $X$ is contained in $[N/4,N/2]$
we see that $x, x+d, x+2d$ is an arithmetic progression in $\Z$ and not just in $\Z_N$.  Theorem \ref{ap-theorem} follows.
\endproof

\begin{remark}
We have in fact shown that the number of triples $(p_1,p_2,p_3)$ of Chen primes with $p_1 + p_3 = 2p_2$ and with each $p_i \leq N$ is $\gg N^2 / \log^6 N$; it is possible to show from \eqref{eq3} that this bound is sharp up to an absolute multiplicative constant. It is clear that if one had a lower bound $\pi_2(N) \gg N /\log^2 N$ for the number of twin primes less than $N$ (this, of course, is one of the most famous open conjectures in prime number theory) then a simple adaptation of our argument would produce infinitely many
triples $(p_1,p_2,p_3)$ of twin primes in arithmetic progression, and the number of such triples less than $N$
would be $\gg N^2 / \log^6 N$.  Similar remarks hold for any other constellation of primes. \end{remark}

\section{Appendix: Construction of the enveloping sieve}\label{appendix1}

\ni In this appendix we construct the enveloping sieve for Proposition \ref{prop3} and verify its properties.  
Some of the material here is rather standard and will not be used directly in the rest of the paper, except via Proposition \ref{prop3}.\vs

\ni Let us recall the setup. We took a fixed polynomial $F$ and defined the ``locally sieved'' sets
\[ X_q = \{n \in \mathbb{Z}_q : (q, F(n)) = 1\}.\]
We begin by analyzing these locally sieved sets $X_q$. From \eqref{rho-def} we have
\begin{equation}\label{x-rho}
\gamma(q) = \E_{n \in \Z_q} \mathbf{1}_{(q,F(n)) = 1}  = \E( \mathbf{1}_{X_q} ) = \frac{1}{q} |X_q|.
\end{equation}
From the Chinese remainder theorem we see that $|X_{qq'}| = |X_q| |X_{q'}|$ whenever
$(q,q') = 1$, and also that $|X_{p^m}| = p^{m-1} |X_{p}|$ whenever $p^m$ is a prime power.  From
this and \eqref{x-rho} we see that $\gamma$ is a multiplicative function, and indeed
\begin{equation}\label{rho-form}
\gamma(q) = \prod_{p|q} \gamma(p).
\end{equation}

\ni Observe that if $q | M$ then the projection $n \mapsto n \md q$ will map $X_M$ to $X_q$.  Indeed,
this map has uniform fibres:

\begin{lemma}\label{uniform}  If $q | M$, then
\[ \E_{n \in \Z_M, n \equiv m \mdsublem q }\mathbf{1}_{X_M}(n)  = \frac{\gamma(M)}{\gamma(q)} \mathbf{1}_{X_q}(m)\]
for all $m \in \Z_q$. 
\end{lemma}

\begin{proof} Let $r$ be the largest factor of $M$ which is coprime to $q$.  We observe that if $n \md q = m$, then
\[\mathbf{1}_{X_M}(n) = \mathbf{1}_{X_q}(m) \mathbf{1}_{X_r}(n \md r).\]
From the Chinese remainder theorem we thus see that
\[ \E_{n \in \Z_M, n \equiv m \mdsub q }\mathbf{1}_{X_M}(n)  = \mathbf{1}_{X_q}(m) \E( \mathbf{1}_{X_r} )
= \mathbf{1}_{X_q}(m) \gamma(r)\]
and the claim follows from \eqref{rho-form}.
\end{proof}\vs

\ni We now analyze the Fourier behaviour of the $X_q$. For each $q \in \mathbb{N}$ and $a \in \mathbb{Z}$ we make the temporary definition
\begin{equation}\label{t-def} t(a,q) := 
 \E_{n \in X_q} e_q(an) = \frac{1}{q \gamma(q)} \sum_{n \in X_q} e_q(an).
\end{equation}
Note that the denominator is non-zero by \eqref{nondeg-def}. 

\begin{lemma}  If $q \in \Z^+$ and $a \in \Z_q^*$, then $t(a,q) = t(ak,qk)$ for any $k \in \Z^+$.  
\end{lemma}

\begin{proof}  Let $k'$ be the largest factor of $k$ which is coprime to $q$.  Observe that
$n \in X_{qk}$ if and only if $n \md q \in X_q$ and $n \md{k'} \in X_{k'}$.  Since
$e_q(an)$ depends only on $n \md q$, we thus see from the Chinese remainder theorem that
\[ \E_{n \in X_{qk}} e_q(an)  = \E_{n \in X_q} e_q(an) ,\]
and the claim follows.
\end{proof}\vs

\ni As a consequence of this we may define the \emph{normalised Fourier coefficients} $s(a/q)$ by
\begin{equation}\label{s-def} s(a/q) := t(a,q).\end{equation} We now record some useful identities and estimates for $s(a/q)$.

\begin{lemma}\label{s-lemma}  Let $q, q' \in \Z^+$ and $a, a' \in \Z$.  Then the following statements hold.
\begin{itemize}
\item[(i)] $s(0) = 1$.
\item[(ii)] $s$ is periodic modulo 1.  In particular when considering $s(a/q)$, we may think 
of $a$ as an element of $\Z_q$.
\item[(iii)] If $(q,q') = 1$, then
$$s(a/q + a'/q') = s(a/q) s(a'/q').$$
\item[(iv)] If $a \in \Z_q^*$ and $q$ is not square-free, $s(a/q) = 0$.
\item[(v)] If $a \in \Z_q^*$, $q > 1$, and $\gamma(q) = 1$, then $s(a/q) = 0$.
\item[(vi)] If $a \in \Z_q^*$, then $|s(a/q)| \leq \prod_{p | q} (1 - \gamma(p))/\gamma(p)$.
\end{itemize}
\end{lemma}

\begin{proof} The first two claims are immediate from \eqref{s-def}.  To confirm (iii), define $b \md{qq'}$ by the congruences $b \equiv aq' \md{q}$ and $b \equiv a'q \md{q'}$, and $r,r'$ by the relations $qr \equiv 1 \md{q'}$, $q'r' \equiv 1 \md{q}$. Then it is easy to check that
\begin{equation}\label{eq7.11} s(a/q)s(a'/q') = \frac{1}{|X_q||X_{q'}|}\sum_{n \in X_q}\sum_{n' \in X_{q'}} e_{qq'} \left(b(q'r' n + qr n')\right).\end{equation}
Observe that $q'r' n + qr n'$ is congruent to $n \md{q}$ and to $n' \md{q'}$. Hence the right-hand side of \eqref{eq7.11} is indeed just
\[ \frac{1}{|X_{qq'}|}\sum_{m \in X_{qq'}} e_{qq'}(bm) = s(b/qq') = s(a/q + a'/q').\]
To prove (iv), we write $q = p^2 k$ for some prime $p$ and integer $k$, and observe that $X_{p^2 k}$ is periodic
mod $pk$.  Since $\frac{a}{q}$ is not an integer divided by $pk$, we see that the expression in \eqref{s-def} necessarily
vanishes. The fifth claim (v) follows by inspection.  Finally, we prove (vi).  By (iv) we may assume $q$ is square-free.  By (iii) we may assume $q$ is prime.  But then (vi) follows by writing
\[ \sum_{n \in X_q} e_q(an) = \sum_{n \in \Z_q} e_q(an) - \sum_{n \in \Z_q \backslash X_q} e_q(an),\]
and observing that the former sum vanishes since $a \in \Z_q^*$, while the latter sum is at most
$|\Z_q \backslash X_q| = q (1-\gamma(q))$ by \eqref{x-rho}.
\end{proof}\vs

\ni Now we can write $\mathbf{1}_{X_{M}}$ as a Fourier series, with coefficients determined by the $\gamma$ and $s$ functions.
\begin{lemma}\label{id-lemma}  For any $M$ and any $n \in \Z_M$, we have
\begin{equation}\label{xr-ident}
 \mathbf{1}_{X_{M}}(n) = \gamma(M) \sum_{q | M} \sum_{a \in \Z_q^*} s(a/q) e_q(-an).
\end{equation}
\end{lemma}

\begin{proof}  The right-hand side can be rewritten as $\gamma(M) \sum_{b \in \Z_M} s(b/M) e_M(-bn)$.
But from \eqref{t-def} and \eqref{s-def} we have $\gamma(M) s(b/M) = \E_{n \in \Z_M} \mathbf{1}_{X_{M}}(n) e_M(bn)$.  The
claim then follows from the Fourier inversion formula.
\end{proof}\vs

\ni Motivated by the identity \eqref{xr-ident}, let us define the function\footnote{While at first glance this expression appears to be complex-valued, it is in fact real, as one can see from conjugation symmetry, or from more explicitly real formulae for $\alpha_R$ below such as \eqref{lambda-form}.} $\alpha_R: \Z \to \R$ by
\begin{equation}\label{lambda-def}
 \alpha_R(n) :=  \sum_{q \leq R} \sum_{a \in \Z_q^*} s(a/q) e_q(-an).
\end{equation}
Thus $\alpha_R$ is some sort of truncated version of $\mathbf{1}_{X_{R!}} / \gamma(R!)$,
which is normalized to have average value 1 (as can be seen by considering the $q=1$ term of the above sum, which
is the only term which does not oscillate in $n$). 

\begin{lemma}\label{G-stable}  Let $h$ denote the multiplicative function
\begin{equation}\label{h-def}
 h(q) := \mu(q)^2  \prod_{p|q} \frac{1 - \gamma(p)}{\gamma(p)}.
\end{equation}
Then for all $n \in X_{R!}$, we have $\alpha_R(n) = G(R)$, where $G(R) = G(F,R)$ is the quantity
\begin{equation}\label{G-def} G(R) := \sum_{q \leq R} h(q).
\end{equation}
\end{lemma}
\begin{remark} Note in particular that $G(R)$ is positive (the $q=1$ summand is equal to 1).\end{remark}

\begin{proof}  Fix $n \in X_{R!}$. Recall the definition \eqref{lambda-def} of $\alpha_R$. In the first instance we will work with the inner sum $\sum_{a \in \Z_q^*} s(a/q) e_q(-an)$ occurring in that definition. We may restrict $q$ to be square-free since only such $q$ contribute to the sum \eqref{lambda-def} in view of Lemma \ref{s-lemma} (iv). We have
\begin{eqnarray*}
\sum_{a \in \mathbb{Z}_q^*} s(a/q) e_q(-an) & = & \sum_{a \in \mathbb{Z}_q^*} \E_{m \in X_q} e_q(am) e_q(-an) \\ & = & \E_{m \in X_q} \sum_{a \in \mathbb{Z}_q} e_q(a(m-n))
 \mathbf{1}_{(a,q)=1}\\
& = &  \E_{m \in X_q} \sum_{a \in \Z_q} e_q(a(m-n)) \sum_{d | (a,q)} \mu(d) \\
& = &
 \sum_{d|q} \mu(d) \E_{m \in X_q} \sum_{a \in \Z_q: d|a} e_q(a(m-n)) \\
& = & 
 \sum_{d|q} \mu(d) \E_{m \in X_q} \frac{q}{d}\mathbf{1}_{m \equiv n \mdsub{q/d}} \\
& = & \sum_{d|q} \frac{\mu(d) \E_{m \in \Z_q, m \equiv n \mdsub {q/d}} \mathbf{1}_{X_q}(m)}{
\E_{m \in X_q}\mathbf{1}_{X_q}(m)  }. \end{eqnarray*}
By Lemma \ref{uniform}, we have
$$
\frac{\E_{m \in \Z_q, m \equiv n \mdsub {q/d}} \mathbf{1}_{X_q}(m) }{
\E_{m \in \Z_q} \mathbf{1}_{X_q}(m) } = \frac{1}{\gamma(q/d)} \mathbf{1}_{X_{q/d}}( n ).$$
Now we substitute into \eqref{lambda-def}, the definition of $\alpha_R$. Since $q$ is being supposed square-free, $d$ and $q/d$ are coprime. Hence
\begin{equation}\label{lambda-form}
 \alpha_R(n) = \sum_{q \leq R} \mu^2(q)\sum_{d|q: (\frac{q}{d},F(n)) = 1} \frac{\mu(d) \gamma(d)}{\gamma(q)}.
\end{equation}
Since $n \in X_R$, we have $(\frac{q}{d}, F(n)) = 1$ for all $q \leq R$ and $d|q$.  
Furthermore from the multiplicativity of $\gamma$ we have
\[ \sum_{d|q} \frac{\mu(d) \gamma(d)}{\gamma(q)} = \prod_{p|q} \frac{1 - \gamma(p)}{\gamma(p)} = h(q),\]
and the claim follows.
\end{proof}\vs

\ni The function $\alpha_R$ is not quite a candidate for an enveloping sieve, because it can be negative.  To resolve
this problem, we shall simply square (and renormalize) $\alpha_R$, defining
\begin{equation}\label{beta-def}
 \beta_R(n) := \frac{1}{G(R)} |\alpha_R(n)|^2.
\end{equation}
Clearly $\beta_R$ is non-negative, has $R!$ as a period,
and equals
$G(R)$ on $X_{R!}$.  We now give a Fourier representation for $\beta_R$ similar to \eqref{lambda-def}, but
with slightly larger Fourier coefficients.

\begin{proposition}\label{rr-prop}\cite{Ramare-Ruzsa}  We have the representation
\[ \beta_R(n) =  \sum_{q \leq R^2} \sum_{a \in \Z_q^*} w(a/q) e_q(-an),\]
where the coefficients $w(a/q)$ obey the bounds
\begin{equation}\label{w-bounds}
 |w(a/q)| \leq 3^{\varpi(q)} |s(a/q)|
\end{equation}
\textup{(}here $\varpi(q)$ denotes the number of prime factors of $q$\textup{)}.  Also we have $w(0) = 1$.
\end{proposition}

\begin{remark} More precise asymptotics for $w(a/q)$ are available, see \cite{Ramare-Ruzsa} (and see
\cite{Ramare} for an even more precise statement in the case $F(n)=n$).  We will not need these refinements here,
however.\end{remark}

\begin{proof}  We follow the arguments of Ramar\'e and Ruzsa \cite{Ramare-Ruzsa}.
Our starting point is the formula \eqref{lambda-form}, which was derived for all $n$ (not just those
in $X_{R!}$).  Writing $q = dr$, we thus have\footnote{Note that the contribution of the $q$ which are not square free can be added in, since it is easy to show that $\sum_{dr = q} \frac{\mu(d)}{\gamma(r)} \mathbf{1}_{(r,F(n))=1} = 0$ for all such $q$.
Indeed, if $p^2 | q$, then the contributions of the cases $p|d, p|r$ and $p \not | d, p^2 | r$ to this sum will 
cancel each other out.},
\begin{equation}\label{qd}
 \alpha_R(n) = \sum_{d, r: dr \leq R} \frac{\mu(d)}{\gamma(r)} \mathbf{1}_{(r,F(n)) = 1}.
\end{equation}
Inserting \eqref{qd} into \eqref{beta-def}, we obtain
\[ \beta_R(n) = \frac{1}{G(R)} \sum_{\substack{d_1,r_1,d_2,r_2\\ d_1r_1, d_2r_2 \leq R}} \frac{\mu(d_1) \mu(d_2)}{\gamma(r_1) \gamma(r_2)} \mathbf{1}_{([r_1,r_2],F(n)) = 1},\]
where $[r_1,r_2]$ denotes the least common multiple of $r$ and $r'$.  Also note that $\alpha_R$, and hence $\beta_R$,
is periodic with period $R!$.  By Fourier inversion, we then have
\[ \beta_R(n) = \sum_{q | R!} \sum_{a \in \Z_q^*} w(a/q) e_q(-an)\]
where
\[ w(a/q) := \frac{1}{G(R)}\E_{n \in \Z_{R!}}
 \sum_{\substack{d_1,r_1,d_2,r_2 \\ d_1r_1, d_2r_2 \leq R}} \frac{\mu(d_1) \mu(d_2)}{\gamma(r_1) \gamma(r_2)} 
\mathbf{1}_{([r_1,r_2],F(n)) = 1} e_q(an).\]
Fix $a, q$.
Observe that the quantity $\mathbf{1}_{([r_1,r_2],F(n)) = 1}$, as a function of $n$, is periodic with period $[r_1,r_2]$.  Thus
$\E_{n \in \Z_{R!}} \mathbf{1}_{([r_1,r_2],F(n)) = 1} e_q(an) $ vanishes unless $q$ divides $[r_1,r_2]$, in which
case the expression is equal to $s(a/q) \gamma([r_1,r_2])$ by \eqref{s-def}.  In particular, since
$[r_1,r_2] \leq R^2$, we know that $w(a/q)$ vanishes when $q > R^2$, and in the remaining cases $q \leq R^2$ we have
\[ w(a/q) = \frac{s(a/q)}{G(R)}
\sum_{\substack{d_1,r_1,d_2,r_2\\ d_1 r_1, d_2 r_2 \leq R\\ q | [r_1,r_2]}} \frac{\mu(d_1) \mu(d_2)}{\gamma(r_1) \gamma(r_2)}
\gamma([r_1,r_2]).\]
From \eqref{rho-form} and \eqref{h-def} we have
$$ \frac{\gamma([r_1,r_2])}{\gamma(r_1) \gamma(r_2)} = \prod_{p | (r_1,r_2)} \frac{1}{\gamma(p)}
= \prod_{p | (r_1,r_2)} (1 + h(p)) = \sum_{t | (r_1,r_2)} h(t),$$
and thus 
\begin{equation}\label{w-form}
 w(a/q) = \frac{s(a/q)}{G(R)} \sum_{t \leq R} h(t)
\sum_{\substack{d_1,r_1,d_2,r_2\\ d_1r_1, d_2r_2 \leq R\\ q | [r_1,r_2]\\ t | (r_1,r_2)}} \mu(d) \mu(d').
\end{equation}
Thus to prove \eqref{w-bounds}, it suffices by \eqref{G-def} to show that
\[
 |\sum_{\substack{d_1,r_1,d_2,r_2\\ d_1r_1, d_2r_2 \leq R\\ q | [r_1,r_2]\\ t | (r_1,r_2)}} \mu(d_1) \mu(d_2)| \leq 3^{\varpi(q)}
\]
for all $t \leq R$.  While this is in fact true for all $q$, it in fact suffices in light of Proposition \ref{s-lemma} (iv) to verify it for square-free $q$, in which case the formulae are slightly simpler.\vs

\ni Fix $t$.  We take advantage of the constraint $t | (r_1,r_2)$ to write $r_j = t \tilde r_j$, $l_j = d_j \tilde r_j$, and $\tilde q := q/(q,t)$, and observe that the constraint $q | [r_1,r_2]$ is equivalent to $\tilde q | [\tilde r_1, \tilde r_2]$.  Also, $\tilde q$ is square-free.  Thus we reduce to proving that
\begin{equation}\label{divisors-bound}
 |\sum_{l_1,l_2 \leq R/t} A(l_1,l_2)| \leq 3^{\varpi(q)}
\end{equation}
where
\[ A(l_1,l_2) := \sum_{\substack{\tilde r_1 | l_1, \tilde r_2 | l_2\\ \tilde q | [\tilde r_1,\tilde r_2]}} \mu(l_1/\tilde r_1) \mu(l_2/\tilde r_2).\]
This expression can be worked out quite explicitly:

\begin{lemma}  Let $\tilde q$ be square-free.  Then the quantity $A(l_1,l_2)$ vanishes 
unless $l_1, l_2 | \tilde q$ and $\tilde q | l_1 l_2$, in which case $|A(l_1,l_2)| = 1$.
\end{lemma}

\begin{proof} It suffices to prove this claim in the ``local'' case when $l_1, l_2, \tilde q$ are all powers of a single prime $p$, since
the general case then follows by splitting $l_1, l_2, \tilde q$ into prime factors and exploiting the multiplicativity of
the summation in $A(l_1,l_2)$.  If $\tilde q=1$ then the constraint $\tilde q | [\tilde r_1,\tilde r_2]$ is vacuous
and the claim follows from M\"obius inversion.  Since $\tilde q$ is square-free, the only remaining case is when
$\tilde q = p$.  In this case, $\tilde q | [\tilde r_1,\tilde r_2]$ if and only if not both of $\tilde r_1$, $\tilde r_2$ equal one, 
and so 
\[ A(l_1,l_2) = \sum_{\tilde r_1 | l_1, \tilde r_2 | l_2} \mu(l_1/\tilde r_1) \mu(l_2/\tilde r_2)
- \mu(l_1) \mu(l_2).\]
Applying M\"obius inversion (or direct computation) we thus see that the expression $A(l_1,l_2)$ equals $-1$ when $(l_1,l_2) = (p,p)$, equals $+1$ when $(l_1,l_2) = (p,1), (1,p)$, and vanishes otherwise, and the claim follows.
\end{proof}\vs

\ni From this lemma we see that there are at most $3^{\varpi(\tilde q)} \leq 3^{\varpi(q)}$ pairs $(l_1,l_2)$
for which $A(l_1,l_2)$ is non-zero, and in each of these cases we have $|A(l_1,l_2)| = 1$.  The claim \eqref{divisors-bound} follows.  Using \eqref{w-form}, this proves \eqref{w-bounds}.\vs

\ni In the case $a=0$ and $q=1$, a closer inspection of the above argument shows that the quantity in absolute values in
\eqref{divisors-bound} is in fact exactly equal to 1.  In view of \eqref{w-form}, Lemma \ref{s-lemma} (i) and \eqref{G-def}, this gives $w(0) = 1$ as desired. This concludes the proof of Proposition \ref{rr-prop}.
\end{proof}\vs

\ni Finally, we rewrite $\alpha_R$ (and hence $\beta_R$) in the more familiar notation of the Selberg sieve. This will allow us to conclude 
the crude estimate for $\beta_R$ demanded in Proposition \ref{prop3} (ii).

\begin{lemma}\label{lem7.8}  We have
\[ \alpha_R(n) = G(R)\sum_{\substack{d \leq R \\ d | F(n)}} \lambda_d,\]
where $\lambda_d$ is defined for squarefree $d$ by
\begin{equation}\label{lam-k-def} \lambda_d := \frac{\mu(d)G_d(R/d)}{\gamma(d)G(R)},\end{equation}
and the quantity $G_d$ is defined by
\[ G_d(x) := \sum_{\substack{q \leq x \\ (q,d) = 1}} h(q).\]
\end{lemma}
\begin{remark} The notation here is reasonably consistent with what is normal in sieve theory, as may be found for example in \cite{Halberstam-Richert}.  For some further discussion of the weights $\lambda_d$ in the special case $F(n)=n$, see Appendix \ref{appendix2} below.\end{remark}

\begin{proof}
We start with \eqref{qd} and use M\"obius inversion to write
\[ \mathbf{1}_{(r,F(n))=1} = \sum_{m | r, m | F(n)} \mu(m),\]
and hence
\[ \alpha_R(n) = \sum_{m \leq R: m | F(n)} \mu(m) \sum_{s, r: sr \leq R, m | r} \frac{\mu(s)}{\gamma(r)}.\]
Writing $r =: mr'$ and $q := sr'$, we thus have
\[ \alpha_R(n) = \sum_{m \leq R: m | F(n)} \mu(m) \sum_{q \leq R/m} \sum_{r' | q} \frac{\mu(q/r')}{\gamma(mr')}.\]
From \eqref{rho-form} a calculation gives
\[ \sum_{r' | q} \frac{\mu(q/r')}{\gamma(mr')} = \frac{\mathbf{1}_{(q,m) = 1}}{\gamma(m)} \prod_{p|q} \mu^2(q)\big(\frac{1}{\gamma(q)} - 1\big)
= \frac{\mathbf{1}_{(q,m) = 1}}{\gamma(m)} h(q),\]
and the claim follows.
\end{proof}

\begin{lemma}\label{lem7.9}
Suppose that $d$ is squarefree. Then $|\lambda_d| \leq 1$.
\end{lemma}
\begin{proof} This is a well-known fact in sieve theory. The proof (which may be found in \cite[p. 100]{Halberstam-Richert}) consists in observing that
\[ G(R) = \sum_{l | d} h(l)G_d(R/l) \geq \big(\sum_{l|d} h(l)\big) G_d(R/d) = \frac{G_d(R/d)}{\gamma(d)},\] which implies the result immediately.\end{proof}\vs

\noindent\textit{Proof of Proposition \ref{prop3}.}  
To prove (i), it suffices in view of Lemma \ref{G-stable} to show that
\[ G(R) \gg_k \frac{ \log^k R}{\mathfrak{S}_F}.\]
In \cite[Lemma 4.1]{Halberstam-Richert}, the bound
\[ G(R) \; \gg_k \; \prod_{p < R} \frac{1}{\gamma(p)}\]
 is obtained. The claim (i) now follows from \eqref{gamma-def} and the classical bound 
$\prod_{k < p \leq R} (1 - k/p) \; \ll_k \; (\log R)^{-k}$. \vs

\ni Let us now consider (ii), that is the statement $\beta_R(n) \ll_{\epsilon} N^{\epsilon}$. By Lemmas \ref{lem7.8} and \ref{lem7.9} we have the bound $\beta(n) \leq G(R)d(F(n))^2$. Recall that $1 - \gamma(p) \leq k/p$ for all $p$ (cf. \eqref{rho-explicit}), together with the non-degeneracy assumption \eqref{nondeg-def}, which implies that $\gamma(p) \geq 1/p$. For any $\epsilon > 0$, we may therefore infer the bound
\begin{equation}\label{eq467} h(q) = \mu(q)^2 \prod_{p|q} \frac{1- \gamma(p)}{\gamma(p)} \leq k^k\mu(q)^2 \prod_{\substack{p|q \\ p > k}} \frac{k/p}{1-k/p} \ll_k k^{\varpi(q)} \log^k q / q
\ll_{k,\epsilon} q^{\epsilon-1}\end{equation}
for any $\epsilon > 0$. Here, we have invoked well-known inequalities for the divisor function (recall that $|a_i|, |b_i| \leq N$, whence $|F(n)| \ll_k N^{2k}$). Using \eqref{eq467}, and recalling \eqref{G-def}, we infer that
\[ G(R) \ll_{k,\epsilon} R^{\epsilon/2} \leq N^{\epsilon/2},\] concluding the proof of (ii). We remark that more precise bounds on $G(R)$ may be found in \cite{Halberstam-Richert}.\vs

\ni To prove (iii) and (iv), it suffices by Lemma \ref{s-lemma} (vi) and Proposition \ref{rr-prop} to show that
\[ 3^{\varpi(q)} \prod_{p|q} (1 - \gamma(p))/\gamma(p) \ll_\epsilon q^{\epsilon-1}\]    
for all square-free $q$. But the left-hand side is just $3^{\varpi(q)}h(q)$, whence by \eqref{eq467} and the classical bound $\varpi(q) \ll \log q /  \log \log q$ we are done.
\endproof

\section{Appendix: A comparison of two enveloping sieves.} \label{appendix2}

\ni Let us begin with a few further remarks concerning what we mean by an enveloping sieve. The idea of considering the sieve $\beta_R$ not simply as a tool for obtaining upper bounds but as an interesting function in its own right is one we got by reading papers of Ramar\'e \cite{Ramare} and Ramar\'e-Ruzsa \cite{Ramare-Ruzsa}. Ramar\'e generously acknowledges that the essential idea goes back to Hooley \cite{Hooley}. It is in those papers (and in fact so far as we know \textit{only} in those papers) that one finds the term \textit{enveloping sieve}.\vs

\ni The name, though perhaps nonstandard, seems to us to be appropriate for describing a situation where one uses a single majorant for the primes to deal simultaneously with a whole family of problems which might normally, by themselves, be dealt with by a family of different sieves.\vs

\ni In Proposition \ref{prop3} take $F(n) = n$, so that $G(R)^{-1}\beta_R$, as defined in \eqref{beta-def}, is a majorant for the primes themselves (or at least the primes $> R$). The expansion of $\beta_R$ as a Fourier series in Proposition \ref{prop3} can be thought of as a way of using a sieve for the primes to tell us about the behaviour of the primes restricted to arithmetic progressions $a \md{q}$. This might normally be handled by a separate procedure for each $(a,q)$.\vs

\ni The authors \cite{green-tao-primes} used a somewhat different majorant for the primes, which we call $\beta'_R$ below (a closely related object was called $\nu$ in \cite{green-tao-primes}). We were concerned with the estimation of correlations

\begin{equation}\label{correlations} \E_{n \leq N} \beta'_R(n+h_1) \dots \beta'_R(n + h_m).\end{equation}

\ni Here, then, a single sieve $\beta'_R$ is being used to deal with configurations $(n+h_1,\dots,n+h_m)$ which might more normally be treated using a Selberg sieve adapted to the polynomial $F(n) = (n + h_1)\dots (n+h_m)$. This, we think, should also qualify as a use of an enveloping sieve.\vs

\ni Let us return to the function $\beta_R(n)$. By Lemma \ref{lem7.8} it is given by the formula

\[ \beta_R(n) = G(R)(\sum_{\substack{d | n \\ d \leq R}} \lamsel_d)^2.\]
$G(R)$ is as defined in \eqref{G-def}. In this particular case we have from \eqref{h-def} that $h(m) = 1/\phi(m)$ and so
\[ G(R) = \sum_{m \leq R} \frac{\mu^2(m)}{\phi(m)}.\] The weights $\lamsel_d$ (SEL stands for Selberg) are defined by \eqref{lam-k-def}; when $F(n) = n$ we have $\gamma(d) = \phi(d)/d$, and thus
\begin{equation}\label{lamsel-def}
 \lamsel_d := \frac{\mu(d) G_d(R/d)}{\gamma(d)G(R)} =
\frac{d \mu(d)}{\phi(d) G(R)}\sum_{\substack{q \leq R/d \\ (q,d) = 1}} 1/\phi(q).
\end{equation}

\ni Why did we consider $\beta_R$? For us, the motivating factor was that the Fourier transform of $\beta_R$ could be described quite accurately, as was done for example in Proposition \ref{rr-prop}. There are many other uses for $\beta_R$, one of the most famous being to the so-called Brun-Titchmarsh problem of proving an upper bound for the number of primes in an interval $[x,x+y)$. The function $\beta_R$ and the associated weights $\lamsel_d$ were first considered by Selberg and constitute the simplest instance of Selberg's \textit{upper bound sieve} (also known as the $\Lambda^2$-sieve).\vs

\ni As is well-known, the weights $\lamsel_d$ arise from a certain natural quadratic optimization problem. Indeed for \emph{any} real weights $\lambda_d$ such that $\lambda_1 = 1$ the function
\[ \psi(n) := (\sum_{\substack{d | n \\ d \leq R}} \lambda_d )^2\] forms a pointwise majorant for the primes $> R$. For a given integer $N$ we might choose $R$ to be a small power of $N$ and then attempt to choose the $\lambda_d$ so that $\sum_{n \leq N} \psi(n)$ is as small as possible. After changing the order of summation one finds that
\begin{equation}\label{eq8.69} \sum_{n \leq N} \psi(n) = N \sum_{d_1,d_2 \leq R} \frac{\lambda_{d_1}\lambda_{d_2}}{[d_1,d_2]} + O(R^2\sup_{d \leq R}\lambda^2_d).\end{equation}
In practice the error will be small, and we become interested in minimising the quadratic form
\[ Q(\mathbf{\lambda}) := \sum_{d_1,d_2 \leq R} \frac{\lambda_{d_1}\lambda_{d_2}}{[d_1,d_2]} = \sum_{\delta \leq R} \phi(\delta) ( \sum_{\delta | d} \frac{\lambda_d}{d} )^2\]
 subject to $\lambda_1 = 1$. $Q(\mathbf{\lambda})$ may be written as $\sum_{\delta \leq R} \phi(\delta) u_{\delta}^2$ where \begin{equation}\label{u-transform} u_{\delta} := \sum_{\delta | d} \lambda_d/d;\end{equation} in these new coordinates the condition $\lambda_1 = 1$ translates to 
 \[ \sum_{\delta \leq R} \mu(\delta) u_{\delta} = 1. \]
Observe that under these conditions we have
\begin{equation}\label{key-eq} Q(\mathbf{\lambda}) = \sum_{\delta \leq R} \phi(\delta) \left(u_{\delta} - \frac{\mu(\delta)}{G(R)\phi(\delta)}\right)^2 + G(R)^{-1},\end{equation}
which makes it completely clear that the minimum of $Q(\lambda)$ is $1/G(R)$ and that this occurs when $u_{\delta} = \mu(\delta)/\phi(\delta)G(R)$. This corresponds to the choice $\lambda_d = \lamsel_d$.\vs

\ni Now it is well-known that $G(R) \sim \log R$. In fact one has the asymptotic
\begin{equation}\label{eq8.70} G(R) = \log R + \gamma + \sum_p \frac{\log p}{p(p-1)} + o(1).\end{equation}
This may be found in \cite{Ramare}, for example; Montgomery \cite{Montgomery2} traces the result back at least as far as Ward \cite{Ward}.\vs

\ni We remark that all this together with Lemma \ref{lem7.9} (which, in our new notation, tells us that $|\lamsel_d| \leq 1$) allows us to recover the fact that if $R = N^{1/2 - \epsilon}$ then
\[ \sum_{n \leq N} \beta_R(n) = N(1 + o_{\epsilon}(1)).\] This also follows quickly from Proposition \ref{prop3} (cf. Lemma \ref{beta-L1}).\vs

\ni As remarked earlier, we had occasion in our paper \cite[Chapter 8]{green-tao-primes} to consider a somewhat different majorant for the primes $> R$. To enable a comparison with $\beta_R$, we write this majorant in the form
\begin{equation}\label{betadash-def} \beta'_R(n) = G(R) (\sum_{\substack{d | n \\ d \leq R}} \lamgy_d)^2\end{equation} where 
\[ \lamgy_d = \frac{\mu(d)\log(R/d)}{\log R}.\]Objects of this type were extensively analysed by Goldston and Y{\i}ld{\i}r{\i}m \cite{gy,gy2,gy3}, and in particular they saw how to asymptotically evaluate certain \textit{correlations} of the form \eqref{correlations} provided $R \leq N^{c_m}$. This was a crucial ingredient in our work in \cite{green-tao-primes}, but we should remark that other aspects of the function $\beta'_R$ were investigated much earlier, and indeed $\beta'_R$ was known to Selberg. It has also featured in works of Friedlander-Goldston \cite{friedlander-goldston}, Goldston \cite{goldston} and Hooley \cite{hooley2}, among others. The consideration of the weights $\lamgy$ and the associated majorant $\beta'_R$ may be motivated by looking at an asymptotic for the weights $\lamsel_d$. In \cite[Lemma 3.4]{Ramare}, for example, one finds the result
\begin{equation}\label{G-funct-asymptotic} G_d(R/d) = \log(R/d) + \gamma + \sum_{p \geq 2} \frac{\log p}{p(p-1)} + \sum_{p | d} \frac{\log p}{p},\end{equation} which leads in view of the definition \eqref{lamsel-def} and \eqref{eq8.70} to the asymptotic
\[ \lamsel_d \sim \frac{\mu(d)\log(R/d)}{\log R} = \lamgy_d.\]
At the very least, this suggests looking at $\beta'_R$ as defined in \eqref{betadash-def} above. The oscillatory effect of the M\"obius function means that there is no \textit{a priori} guarantee that $\beta'_R$ is a close approximation to $\beta_R$. It turns out however that $\beta_R$ and $\beta'_R$ are, for many purposes, rather similar.\vs

\ni One way to see this is to inspect the quadratic form $Q(\lamgy)$. The asymptotic
\[ Q(\lamgy) = \frac{1}{\log R} + O(\frac{1}{\log^2 R})\] follows from a result of Graham \cite{graham} and can easily be proved using Goldston and Y{\i}ld{\i}r{\i}m's method (as exposited, for example, in \cite[Ch. 9]{green-tao-primes}).  Comparing with \eqref{key-eq} and using the asymptotic \eqref{eq8.70} one sees that
\[ \sum_{\delta \leq R} \phi(\delta)(u^{\mbox{\scriptsize GY}}_{\delta} - u^{\mbox{\scriptsize SEL}}_{\delta})^2 = O(\frac{1}{\log^2 R}).\] Here, of course, $u^{\mbox{\scriptsize GY}}$ and $u^{\mbox{\scriptsize SEL}}$ are related to $\lamgy$ and $\lamsel$ respectively by the transformation \eqref{u-transform}. If $R \ll N^{1/2 - \epsilon}$ then this translates via \eqref{eq8.69} to a bound
\[ \sum_{n \leq N} \big( \sum_{\substack{d | n \\ d \leq R}} \lamgy_d - \lamsel_d\big)^2 = O(\frac{N}{\log^2 R}),\] which implies that
\[ \sum_{n \leq N} \left(\beta^{\prime 1/2}_R(n) - \beta_R^{1/2}(n)\right)^{2} = O(\frac{N}{\log R}).\]
Using the Cauchy-Schwarz inequality we now get
\begin{eqnarray*}
\Vert \beta'_R - \beta_R \Vert_1 & := & \E_{1 \leq n \leq N} | \beta'_R(n) - \beta_R(n)|  \\ & \leq & \big(\E_{1 \leq n \leq N} (\beta_R^{\prime 1/2}(n) - \beta_R^{1/2}(n))^2 \big)^{1/2} \times \\ & & \qquad \times \;
\big(\E_{1 \leq n \leq N}(\beta_R^{\prime 1/2}(n) + \beta_R^{1/2}(n))^2 \big)^{1/2} \\ & = & O(\frac{1}{\log R})^{1/2} \cdot (\Vert \beta_R\Vert_1 + \Vert \beta'_R\Vert_1)^{1/2} \\ & = & O(\frac{1}{\sqrt{\log R}}).
\end{eqnarray*}
On the other hand, from Lemma \ref{beta-L1} we have $\| \beta_R \|_1 = O(1)$.  From the triangle inequality we thus obtain

\begin{lemma}
Suppose that $R \ll N^{1/2 - \epsilon}$. Then $\beta_R$ and $\beta'_R$ are close in that they satisfy the $L^1$-estimate $\Vert \beta'_R - \beta_R \Vert_1 = O((\log R)^{-1/2})$. \endproof
\end{lemma}

\ni Such an estimate is not, in itself, enough to allow us to deduce enough information about the Fourier expansion of $\beta'_R$ to prove a restriction theorem analogous to that obtained in \S \ref{majorant-sec} for $\beta_R$. Such Fourier information could, we believe, be obtained by an analytic method similar to that in \cite{gy3}, but there seems little use for such a result. Going in the other direction, one might ask whether it is possible to obtain combinatorial information concerning correlations of the form \eqref{correlations} for $\beta_R$. In view of \eqref{G-funct-asymptotic} it looks as though any such effort might require one to address the function
\[ H_R(n) :=  \sum_{\substack{d | n \\ d \leq R}} \mu(d),\] concerning oneself in particular with the $2m^{\mbox{\scriptsize th}}$ moment of this function. When $m = 1$ this was investigated by Dress, Iwaniec and Tenenbaum \cite{DIT} and when $m = 2$ by Motohashi \cite{motohashi}. At this point there are already some rather thorny issues involved, and it does not seem worth the effort of pursuing the matter simply to obtain correlation estimates for $\beta_R$. The function $H_R$ is of some interest in its own right, however, and had an auxilliary role for instance in the ground-breaking work \cite{FI}.\vs

\ni Let us conclude with the following remark. As far as enveloping sieves for the primes are concerned, it seems that both $\beta_R$ and $\beta'_R$ have their own role to play. $\beta_R$ is useful if one wants to do harmonic analysis, whereas $\beta'_R$ is far more appropriate if one wishes to do combinatorics. The two majorants are, however, rather closely related.

\end{document}